\documentclass[ijoc,sglanonrev]{informs4}
\usepackage{eqndefns-left} 
\RequirePackage{tgtermes}
\RequirePackage{newtxtext}
\RequirePackage{newtxmath}
\RequirePackage{bm}

\OneAndAHalfSpacedXII 

\allowdisplaybreaks[2] 

\makeatletter
\AtBeginDocument{%
  \let\INFORMSthebibliography\thebibliography
  \renewcommand{\thebibliography}[1]{%
    \INFORMSthebibliography{#1}%
    \fontsize{9}{10}\selectfont
    \setlength{\itemsep}{0pt}}}
\makeatother
\newlength{\modelskip}
\newlength{\modelrowsep}
\newenvironment{model}
  {\addtolength{\abovedisplayskip}{\modelskip}%
   \addtolength{\belowdisplayskip}{\modelskip}%
   \addtolength{\abovedisplayshortskip}{\modelskip}%
   \addtolength{\belowdisplayshortskip}{\modelskip}%
   \subequations}
  {\endsubequations}

\usepackage{algorithm}
\usepackage{algpseudocode}
\usepackage{tikz}
\usepackage{amsmath,amssymb,amsfonts}
\usepackage[bookmarks=true,linkcolor=red,citecolor=blue,urlcolor=blue,colorlinks=true, breaklinks]{hyperref}
\usepackage{nicematrix}
\usepackage{booktabs}
\usepackage{array}
\usepackage{capt-of}
\usepackage{threeparttable}
\usepackage{pgfplots}
\usepackage{bookmark}
\usepackage{multirow}
\usepackage[normalem]{ulem}
\usepackage{longtable}
\usepackage{subcaption}
\usepackage{lscape}
\usepackage{makecell}
\usepackage{rotating}
\usepackage{fancyvrb}
\usepackage{anyfontsize} 
\usetikzlibrary{positioning,arrows.meta,shapes.geometric,fit,calc}
\usepackage{comment}
\usepackage{xcolor}

\newif\ifREVMARK
\REVMARKfalse
\newcommand{\Fabio}[1]{\ifREVMARK\textcolor{blue}{#1}\else#1\fi}
\newcommand{\red}[1]{#1}
\newcommand{\blue}[1]{\ifREVMARK\textcolor{blue}{#1}\else#1\fi}
\newcommand{\orange}[1]{\ifREVMARK\textcolor{orange}{#1}\else#1\fi}
\newcommand{\magenta}[1]{\ifREVMARK\textcolor{magenta}{#1}\else#1\fi}
\newcommand{\algname}[1]{\texttt{#1}}
\newenvironment{myproof}{
\noindent\textit{Proof.} }{\hfill $\square$}

\newcounter{examplecounter}
\renewcommand{\theexamplecounter}{\arabic{examplecounter}}

\usepackage{natbib}
 \bibpunct[, ]{(}{)}{,}{a}{}{,}%
 \def\bibfont{\small}%
\EquationsNumberedThrough    

\TheoremsNumberedThrough     
\ECRepeatTheorems  %

\MANUSCRIPTNO{}

\makeatletter

\renewcommand{\theARTICLETOP}{}

\renewcommand{\theARTICLEABSTRACT}{%
  \HOOKb
  \vspace*{18pt}%
  \noindent
  \begin{minipage}[t]{\textwidth}
    \parindent1em
    \ABSfont

    \noindent
    \theABSTRACT
    \endgraf

    \vskip5pt

    \theFUNDING
    \theKEYWORDS
    \theSUBJECTCLASS
    \theAREAOFREVIEW
    \theMSCCLASS
    \theORMSCLASS

    \if@BLINDREV
    \else
      \theHISTORY
    \fi

    \noindent\hrulefill
  \end{minipage}%
  \vspace*{0pt}%
}

\renewcommand{\setoddRH}{%
  \hbox to \textwidth{%
    \fs.7.8.%
    \tabcolsep0pt
    \begin{tabular*}{\textwidth}[b]{%
      @{}l@{\extracolsep{\fill}}r@{}%
    }
      {\theRRHFirstLine}
      &
      {\fs.10.10.\thepage}
      \\[-4pt]
      \multicolumn{2}{@{}l@{}}{%
        \VRHDW{0.5pt}{0pt}{\textwidth}%
      }
    \end{tabular*}%
  }%
}

\renewcommand{\setevenRH}{%
  \hbox to \textwidth{%
    \fs.7.8.%
    \tabcolsep0pt
    \begin{tabular*}{\textwidth}[b]{%
      @{}l@{\extracolsep{\fill}}r@{}%
    }
      {\fs.10.10.\thepage}
      &
      {\theLRHFirstLine}
      \\[-4pt]
      \multicolumn{2}{@{}l@{}}{%
        \VRHDW{0.5pt}{0pt}{\textwidth}%
      }
    \end{tabular*}%
  }%
}

\makeatother

\begin{document}


\RUNAUTHOR{Zhang, Chen, Furini, and Ljubić}


\RUNTITLE{An Exact Algorithm for the Max--Min Covering Location Blocker Problem}


\TITLE{An Exact Algorithm for the Max--Min Covering Location Blocker Problem}

\ARTICLEAUTHORS{%
\AUTHOR{Yun-Tian Zhang\textsuperscript{a}, Chen Chen\textsuperscript{b}, Fabio Furini\textsuperscript{c}, Ivana Ljubić\textsuperscript{d}}
\AFF{\textsuperscript{a}School of Automation, Beijing Institute of Technology, \EMAIL{zhangyuntian@bit.edu.cn}\\ \textsuperscript{b}School of AI, Beijing Institute of Technology, \EMAIL{xiaofan@bit.edu.cn}\\ \textsuperscript{c}Department of Computer, Control and Management Engineering ``Antonio Ruberti'', Sapienza University of Rome, \EMAIL{fabio.furini@uniroma1.it}\\ \textsuperscript{d}Institute for Statistics and Mathematics, Vienna University of Economics and Business, \EMAIL{ivana.ljubic@wu.ac.at}}
} 

\ABSTRACT{%
We introduce the Max–Min Covering Location Blocker Problem, a bilevel optimization problem in which a leader blocks a minimum-cost set of candidate locations so that the optimal coverage of a budget-constrained follower does not exceed a prescribed target. Each customer’s coverage is determined by the least favorable open facility that can serve it. The induced follower coverage set function is nonmonotone and nonsubmodular, so the validity of interdiction cuts cannot be inherited from existing frameworks. We develop an exact nested decomposition that exploits the max–min coverage structure at both levels. At the outer level, we establish the validity of interdiction cuts from this structure and strengthen them with facility-specific coefficients bounding the coverage lost when a critical facility is removed. We solve the NP-hard follower problem by branch-and-Benders-cut after projecting out the customer coverage variables, and characterize an optimal dual solution of the separation problem in closed form, generating Benders optimality cuts in linear time without solving a linear program. Computational experiments on benchmark instances show that both ingredients substantially improve the performance of the method. With a moderate number of candidate locations, instances with 300,000 customers are solved to optimality within one hour. A comparison with two general-purpose bilevel solvers from the literature shows the proposed method to be one to three orders of magnitude faster, and to solve to optimality instances that neither of them closes. The model applies wherever the planner cannot control which facility serves a customer, as in public automated external defibrillator (AED)   networks. A case study on a    network using real AED data from Virginia Beach illustrates the model.
}%
\FUNDING{This research was supported by the National Natural Science Foundation of China under Grant 62673059 and has been conducted during research stay of Y.T. Zhang at La Sapienza, Rome. }



\KEYWORDS{Combinatorial optimization, bilevel problem, covering location problem, blocker problem, Benders decomposition} 

\maketitle



\section{Introduction}
Covering models are widely used to design systems in which facilities provide spatial or service access to geographically distributed customers. In many applications, however, nominal coverage is only part of the planning problem: facilities may become unavailable because of failures, disruptions, accessibility restrictions, or deliberate interdiction, and the system must retain an acceptable service level after such losses. This motivates blocker and interdiction models that identify a minimum-cost set of facilities whose removal is sufficient to reduce the best attainable service level below a prescribed target.

We study this question for a covering system with a max--min service criterion \red{of  \citet{chan2016optimizing}}. A budget-constrained decision maker selects facilities, but the realized service origin for a customer cannot be assumed to be the most favorable selected facility. Instead, the guaranteed coverage of a customer is the minimum service level among the selected facilities capable of serving it. This \red{max--min coverage} structure captures settings in which access, failover, dispatch, or assignment outcomes are not completely controlled by the planner.

The resulting \emph{Max--Min Covering Location Blocker Problem} (MMCLBP) is a bilevel optimization problem. A leader first blocks a minimum-cost set of candidate locations. The follower then optimizes a budget-constrained facility configuration over the surviving locations to maximize total guaranteed coverage. The leader seeks a blocking decision for which this optimal follower value does not exceed a target. A direct solution is computationally challenging because the follower's NP-hard problem must be solved repeatedly and, in the applications targeted here, the number of customers can be orders of magnitude larger than the number of candidate facilities.

Our solution approach exploits the particular max--min \red{coverage} structure at both levels of a nested exact decomposition, and the computational study shows that the two structural ingredients it yields are precisely the components responsible for the practical performance of the final algorithm.

\subsection{Applications and motivating example}\label{sec:applications}
A primary motivating application is the resilience analysis of automated external defibrillator (AED) networks. Sudden cardiac arrest is a time-critical emergency in which rapid access to an AED can be life-saving. Yet an apparently nearby service option may be temporarily inaccessible, already committed, unavailable, or not selected by the operational process. A lay responder does not know where the AEDs are and must search for one with the guidance of a dispatcher and of whatever signage is available, so the best admissible device cannot be assumed to be the one retrieved. For this reason \citet{chan2016optimizing} bracket single-responder behavior between a best case and a worst case, and credit each demand point, in the worst case, with the service level of the least favorable open facility able to reach it. We adopt that criterion. AED-location studies have considered candidate sets ranging from public-access businesses and points of interest to large building-level sets, illustrating the variety of spatial design settings in which coverage and accessibility must be balanced~\citep{sun2016overcoming,tierney2018novel,pourghaderi2022maximum}.

In an AED-drone setting, candidate emergency-service facilities can also act as drone bases, and historical out-of-hospital cardiac-arrest events define spatial demand~\citep{boutilier2022drone,boutilier2017optimizing}. The blocker model then identifies combinations of bases whose simultaneous unavailability is sufficient to reduce the best achievable guaranteed coverage below a prescribed threshold. Importantly, the facilities that are attractive in a nominal deployment need not coincide with the facilities that are most important for preserving reconfiguration flexibility after disruptions. Section~\ref{sec:case-study} develops this application using real data from Virginia Beach.

The MMCLBP can also represent other covering systems in which the realized service origin is not fully controllable, including emergency services, telecommunications, supply chains, and transportation systems. For example, requests in content-delivery networks may be redirected to less favorable edge nodes when preferred servers are unavailable~\citep{zolfaghari2020content}. The common modeling feature is that system performance depends on a conservative guarantee over admissible service origins rather than on the best available origin.

\subsection{Related literature}\label{sec:literature}
The literature most closely related to this work can be organized into three streams. The first concerns covering location models and their extensions. The classical maximal covering location problem selects a prescribed number of facilities to maximize the demand covered within a given service standard~\citep{church1974maximal}. Subsequent extensions relax the all-or-nothing notion of coverage by accounting for service uncertainty or partial coverage~\citep{church2018probabilistic}. For example,~\citet{daskin1983maximum} introduces expected coverage to account for facility unavailability, while~\citet{karasakal2004maximal,alvarez2019exact} allow coverage to decrease gradually with distance. Closest to our setting, \citet{chan2016optimizing} propose three probabilistic coverage models that differ in the assumed AED-retrieval behavior of lay responders. Their worst-case single-responder model evaluates a deployment through the least favorable open facility able to serve each demand point, and they show that its optimal value lower-bounds those of the other two models and of the maximal covering location problem. The max--min partial-coverage criterion used in this paper is that model, together with the two families of linear constraints they derive for it. These models motivate the use of partial or probabilistic coverage measures, but they retain a single decision-making level. The MMCLBP keeps the criterion of \citet{chan2016optimizing} and embeds it in a bilevel blocker setting.

The second stream concerns facility disruption, interdiction, and bilevel location models~\citep{church2018disruption}. \citet{church2004identifying} introduce median and covering facility interdiction models for identifying facilities whose loss is most damaging to a service system. Related protection models explicitly represent the interaction between defensive and disruptive decisions through bilevel optimization~\citep{scaparra2008bilevel}. In the covering context,~\citet{o2011designing} develop a location--interdiction model that accounts for worst-case facility losses when designing a robust coverage network. Methodologically, the closest works are those that solve interdiction games by single-level reformulations built on interdiction cuts, generated dynamically as in the progressive approximation approach of \citet{contardo2022progressive}. \citet{fischetti2019interdiction} tighten such cuts under a \emph{downward monotonicity} assumption, and \citet{taninmis2022branch} extend them to followers maximizing a monotone submodular set function, testing their branch-and-cut on a weighted maximal covering interdiction game. Our interdiction cuts have the same generic form, but neither framework covers the MMCLBP: the max--min coverage function is neither monotone nor submodular, and downward monotonicity fails precisely because removing a facility from a follower solution can \emph{raise} the guaranteed coverage of a customer. \red{Appendix~\ref{ec:submodularity} gives counterexamples for both the follower objective and the optimal follower value.} This is why the loss coefficients of Section~\ref{sec:bigM} must be derived from the max--min structure rather than from marginal gains. Our problem also differs from the location--interdiction models above in both the decision sequence and the service measure: the leader seeks a minimum-cost blocking set that forces the optimal response of a budget-constrained follower below a prescribed target, while follower performance is determined by a max--min partial-coverage function.

The third stream concerns exact decomposition methods for large-scale covering problems. Large-scale covering models have motivated decomposition approaches that project customer-level variables and exploit Benders-type reformulations; see, for example, the covering-specific decomposition of~\citet{cordeau2019benders}. More broadly, modern single-tree and branch-and-Benders-cut implementations exploit callbacks to integrate Benders separation within branch-and-bound~\citep{ljubic2012exact,fischetti2016benders,fischetti2017redesigning}. Most recently, \citet{fadda2026tailored} develop a tailored branch-and-Benders-cut framework for backup covering problems. These works provide the methodological background for our branch-and-Benders-cut algorithm.

\subsection{Contributions and organization}\label{sec:contributions}
The paper makes three main contributions. First, taking the max--min partial-coverage criterion from \citet{chan2016optimizing}, we formulate the MMCLBP as the resilience counterpart of that covering model and develop an exact \red{nested decomposition} for its bilevel structure. The coverage criterion is therefore adopted from the literature; the bilevel formulation and the \red{developed} solution algorithms are new. Second, we exploit the max--min coverage function to derive facility-specific loss coefficients that strengthen the outer interdiction cuts and substantially reduce the number of follower problems that must be solved. Third, we derive \red{a Benders decomposition for the lower-level problem whose optimality cuts are globally valid and} can be separated in closed form in linear time, without solving an auxiliary LP. The computational study isolates the impact of these ingredients, evaluates scalability with respect to the customer dimension, illustrates the model on the Virginia Beach AED case study, and compares the method with two general-purpose bilevel solvers. \Fabio{It also delimits the regime in which the method is effective, with the facility rather than the customer dimension as the binding limit.}

The remainder of the paper is organized as follows. Section~\ref{sec:Problem} defines the MMCLBP and its follower formulation. Section~\ref{sec:outer-decomposition} develops the \red{outer decomposition} and strengthened interdiction cuts. Section~\ref{sec:Benders} develops the \red{branch-and-Benders-cut algorithm} with closed-form separation for the follower, \red{and combines the two into \algname{OD+B\&BC}, the exact algorithm proposed here, whose name records the outer decomposition and the branch-and-Benders-cut used for the lower level}. Section~\ref{sec:DoE} presents the computational study and the Virginia Beach application. Section~\ref{sec:conclusions} concludes. Appendix~\ref{app:general-purpose} reports the comparison with general-purpose bilevel solvers.

\section{The Max--Min Covering Location Blocker Problem}\label{sec:Problem}
We first define the leader--follower interaction, the max--min coverage measure, and the main notation, and then derive a mixed-integer linear programming (MILP) formulation of the follower problem, which is the starting point for both algorithms developed in the following sections.

\subsection{Problem definition and notation}\label{sec:problem-definition}
Let $I=\{1,2,\dots,m\}$ denote the set of $m$ candidate facility locations, and let $J=\{1,2,\dots,n\}$ denote the set of $n$ customers. Each location $i\in I$ has a nonnegative installation cost $f_i\in\mathbb{R}_{\geq 0}$ and a nonnegative blocker cost $e_i\in\mathbb{R}_{\geq 0}$. Each customer $j\in J$ has a \Fabio{positive} importance weight $d_j\in\mathbb{R}_{> 0}$. Let $p_{ij}\in[0,1]$ denote the coverage level provided by a facility at location $i$ to customer $j$. Depending on the application context, $p_{ij}$ can represent a coverage probability or a normalized service quality level. Let \(I_j=\{i\in I: 0<p_{ij}\le 1\}\) denote the set of locations that can provide positive coverage to customer \(j\in J\). Similarly, let \(J_i=\{j\in J: 0<p_{ij}\le 1\}\) denote the set of customers that can be covered by location \(i\in I\). Let $b\in \mathbb{R}_{\geq 0}$ denote the facility-opening budget. Let $t\in\mathbb{R}_{\geq 0}$ denote the \emph{target coverage} threshold.

We introduce the MMCLBP, which belongs to the class of blocker problems, a special class of bilevel optimization problems~\citep{beck2026linear}. The blocker policy selects locations to remove first; facility locations are then optimized in response. We model this interaction as a static Stackelberg game between a leader and a follower. \red{Throughout the paper we use the terms \emph{upper-level problem} and \emph{leader problem} interchangeably, and likewise \emph{lower-level problem} and \emph{follower problem}.} Let $\boldsymbol{x}\in \{0,1\}^m$ denote the vector of upper-level binary variables, where $x_i=1$ if and only if location $i\in I$ is blocked. A location is said to be blocked if the leader makes it unavailable for subsequent facility location by the follower. We use the term block consistently for this upper-level action throughout the paper. Let $\boldsymbol{y}\in \{0,1\}^m$ denote the vector of lower-level binary variables, where $y_i=1$ if and only if a facility is deployed at location $i\in I$. Let the continuous variable $r_j\geq 0$ denote the overall level of coverage provided to customer $j\in J$. We define $I(\boldsymbol{y})=\{i \in I: y_i=1\}$ as the set of selected facility locations and $I_j(\boldsymbol y)=I_j\cap I(\boldsymbol{y})$. Finally, let $\xi(\boldsymbol{x})$ denote the optimal follower coverage value induced by a blocker decision $\boldsymbol{x}$,  that is, the \emph{optimal-value function} of the lower-level problem.   The MMCLBP is defined by the following upper-level and lower-level problems.
The upper-level problem is
{\addtolength{\abovedisplayskip}{\modelskip}%
 \addtolength{\belowdisplayskip}{\modelskip}%
 \begin{equation}
 \min\left\{~\sum_{i\in I} e_i x_i :~~ \xi(\boldsymbol{x})\leq t,~~ x_i\in\{0,1\},~\ i\in I ~\right\}. \label{OF_obju}
 \end{equation}}
For a fixed blocker decision $\boldsymbol{x}$, the \red{value of $\xi(\boldsymbol{x})$} is obtained from the lower-level problem
\begin{model}
\label{OF_lower}
\begin{align}
&& \xi(\boldsymbol{x}) = \max ~~ \sum_{j\in J} d_j r_j \label{OF_objl} \\[\modelrowsep]
&& {\rm subject~to} \quad r_j &=
\begin{cases}
\displaystyle \min_{i\in I_j(\boldsymbol y)} p_{ij}, & \text{if } I_j(\boldsymbol y)\neq\varnothing,\\
0, & \text{otherwise},
\end{cases}
& j\in J, \label{OF_cons_l1} \\[\modelrowsep]
&& \sum_{i\in I} f_i y_i &\leq b, & \label{OF_cons_l2} \\[\modelrowsep]
&& y_i &\leq 1-x_i, & i\in I, \label{OF_cons_l3} \\[\modelrowsep]
&& y_i &\in \{0,1\}, & i\in I, \label{OF_cons_l4} \\[\modelrowsep]
&& r_j &\geq 0, & j\in J. \label{OF_cons_l5}
\end{align}
\Fabio{In this model, the upper-level~\eqref{OF_obju} seeks a minimum-cost subset of locations to be blocked, such that the maximum total worst-case coverage of demand points~\eqref{OF_objl} does not exceed a given threshold $t$.} Constraint~\eqref{OF_cons_l1} defines \(r_j\) as the minimum level of coverage among the selected facilities that can cover customer \(j\), and sets \(r_j=0\) if no such facility is selected. Constraint~\eqref{OF_cons_l2} is a knapsack-like constraint that ensures the available budget $b$ for facility location is not exceeded. Linking constraints~\eqref{OF_cons_l3} force $y_i=0$ at every blocked location. \orange{Note that~\eqref{OF_cons_l2} bounds the budget from above rather than requiring it to be spent exactly. The follower-feasible set is therefore downward closed: every subset of a feasible configuration is itself feasible.}
\end{model}

\subsection{\red{Lower-level formulation}}\label{sec:follower-formulation}
\red{Formulation~\eqref{OF_lower} is nonlinear: while \eqref{OF_objl} is linear in $\boldsymbol r$, each $r_j$ is defined in~\eqref{OF_cons_l1} as a minimum over the index set $I_j(\boldsymbol y)$, which itself is determined by $\boldsymbol y$.} However, two families of linear constraints~\citep{chan2016optimizing} \red{can be used to} provide an equivalent representation in the lower-level maximization problem. First, the coverage variable $r_j$ cannot exceed the coverage level induced by any selected facility that can serve customer $j$:
\begin{equation}
    r_j\leq 1-(1-p_{ij})y_i,~ j\in J,~ i\in I_{j}  \label{linear_1}
\end{equation}
Second, the following sum-based upper bound forces $r_j$ to zero when no covering facility is selected:
\begin{equation}
    r_j\leq \sum\limits_{i\in I_{j}}p_{ij}y_i,~ j\in J. \label{linear_2}
\end{equation}
Together with the maximization objective and \Fabio{$d_j\in\mathbb{R}_{> 0}$}, constraints~\eqref{linear_1}--\eqref{linear_2} make $r_j$ equal to the minimum coverage level among the deployed facilities serving customer $j$, and set it to zero when no such facility is deployed.

For a fixed blocker decision $\tilde{\boldsymbol x}$, the linearized follower problem is obtained by replacing~\eqref{OF_cons_l1} with~\eqref{linear_1}--\eqref{linear_2} and imposing the linking restrictions $y_i\leq 1-\tilde x_i$, $i\in I$, directly. This is the formulation used throughout the algorithms and computational experiments.  The lower-level problem is NP-hard because it contains the maximum coverage location problem (MCLP) as a special case: set $f_i=1$ for all $i$, $b=q$ for a cardinality bound $q$, and $p_{ij}\in\{0,1\}$; the follower then selects at most $q$ facilities to maximize the total demand covered~\citep{megiddo1983maximum}.

\section{\red{Decomposition of the Upper-Level Problem}}\label{sec:outer-decomposition}
This section develops the outer decomposition, the exact framework used to solve the upper-level problem of the MMCLBP, in which the relaxed leader problem is solved repeatedly and refined at each iteration with information obtained from the follower.

\subsection{\red{Single-level reformulation}}\label{sec:vf-reformulation}
Let $Y=\{\boldsymbol y\in\{0,1\}^m:\sum_{i\in I}f_i y_i\leq b\}$ denote the set of follower-feasible facility-opening decisions before imposing the blocking constraints~\eqref{OF_cons_l3}. For any $\boldsymbol y\in Y$, let $r_j(\boldsymbol y)$ denote the coverage induced for customer $j$ according to~\eqref{OF_cons_l1}, and let $\varphi(\boldsymbol y)=\sum_{j\in J}d_jr_j(\boldsymbol y)$ denote the corresponding total coverage value. In words, $Y$ contains every configuration the follower could afford if nothing were blocked and $\varphi(\boldsymbol y)$ is the total weighted guaranteed coverage of a configuration, so that $\xi(\boldsymbol x)=\max\{\varphi(\boldsymbol y):\boldsymbol y\in Y,\ y_i\leq 1-x_i,\ i\in I\}$. \red{Our single-level reformulation} exploits the fact that each fixed follower solution $\tilde{\boldsymbol y}\in Y$ induces an inequality in the space of the upper-level variables. If an upper-level decision $\boldsymbol{x}$ blocks at least one facility in $I(\tilde{\boldsymbol y})$, then $\tilde{\boldsymbol y}$ is no longer a feasible follower response. \Fabio{For each $i\in I(\tilde{\boldsymbol y})$, let $M_i(\tilde{\boldsymbol y})\geq0$ denote a valid loss coefficient, which may depend both on the facility $i$ and on the reference solution $\tilde{\boldsymbol y}$; two specific choices are developed in Section~\ref{sec:bigM}.} 
The resulting single-level reformulation is
\begin{model}
\label{F}
\begin{align}
&& \min ~~ \sum_{i\in I} e_i x_i \label{F_obj} \\[\modelrowsep]
&& {\rm subject~to} \quad t &\geq \varphi(\tilde{\boldsymbol y})
      -\sum_{i\in I(\tilde{\boldsymbol y})} \Fabio{M_i(\tilde{\boldsymbol y})} x_i,
& \tilde{\boldsymbol y}\in Y, \label{F_cons1} \\[\modelrowsep]
&& x_i &\in \{0,1\}, & i\in I. \label{F_domain}
\end{align}
\end{model}
Inequality~\eqref{F_cons1} can be read as follows: if no facility used by $\tilde{\boldsymbol y}$ is blocked, the follower can still attain the coverage $\varphi(\tilde{\boldsymbol y})$, which must then not exceed $t$; each blocked facility $i$ used by $\tilde{\boldsymbol y}$ relaxes this requirement by \Fabio{$M_i(\tilde{\boldsymbol y})$}, an upper bound on the coverage the follower loses without $i$. \red{Following the interdiction literature, we refer to the inequalities~\eqref{F_cons1} as \emph{interdiction cuts}, see, e.g., Section 10.5 in \citet{beck2026linear}.} Since $Y$ may contain exponentially many binary vectors, the complete \red{single-level reformulation} contains exponentially many such inequalities. Therefore, instead of enumerating all follower-feasible solutions, we generate these inequalities dynamically, as described in Section~\ref{sec:outer-algorithm}.

\red{The equivalence between the complete family~\eqref{F_cons1} and the original bilevel
condition $t\geq\xi(\boldsymbol x)$ does not follow from the general results for interdiction cuts
(see, e.g., Section~10.5 in \citet{beck2026linear}), which require the big-M coefficients to be independent on $\tilde{\boldsymbol y}$. The next proposition establishes it under a condition on the coefficients
$M_i(\tilde{\boldsymbol y})$ that the max--min structure allows us to verify directly.}

\begin{proposition}
\label{prop:validity}
\red{Let the coefficients $M_i(\tilde{\boldsymbol y})$ satisfy
$$
\varphi(\tilde{\boldsymbol y}^{\,S})
\geq
\varphi(\tilde{\boldsymbol y})-\sum_{i\in S}M_i(\tilde{\boldsymbol y}),
\qquad
\tilde{\boldsymbol y}\in Y,\ S\subseteq I(\tilde{\boldsymbol y}),
$$
where $\tilde{\boldsymbol y}^{\,S}$ is obtained from $\tilde{\boldsymbol y}$ by setting
$\tilde y_i=0$ for $i\in S$. Then $\boldsymbol x\in\{0,1\}^m$ satisfies every
inequality~\eqref{F_cons1} if and only if $t\geq\xi(\boldsymbol x)$.}
\end{proposition}

\red{The proof is given in Section~\ref{ec:proofs} of the electronic companion.}

\magenta{Although the coverage function $\varphi(\boldsymbol y)$ is nonmonotone
in the selected facilities, the optimal follower value
$\xi(\boldsymbol x)$ cannot increase when additional facilities are
blocked: blocking only restricts the follower's feasible set, and an
available facility need not be selected. The validity of
\eqref{F_cons1} relies on bounding the loss from deleting blocked
facilities from a reference solution, which remains feasible under
the upper-budget constraint. Any coverage increase caused by these
deletions only makes the resulting lower bound more conservative.}

\subsection{Strengthened solution-dependent coefficients}\label{sec:bigM}
Since $0\leq r_j(\boldsymbol y)\leq 1$ for every $j\in J$, the follower value satisfies $\varphi(\boldsymbol y)\leq\sum_{j\in J}d_j$.
 We first derive a simple uniform coefficient that yields a no-good interdiction cut and then exploit the max--min structure to obtain stronger facility-specific coefficients. \Fabio{In particular, taking $M_i(\tilde{\boldsymbol y}):=\varphi(\tilde{\boldsymbol y})$ for every $i\in I(\tilde{\boldsymbol y})$ makes~\eqref{F_cons1} the valid family of interdiction cuts}
\begin{equation}
    t \geq 
    \varphi(\tilde{\boldsymbol y})
    -
    \Fabio{\varphi(\tilde{\boldsymbol y})}
    \sum_{i\in I(\tilde{\boldsymbol y})}x_i,
    \qquad \tilde{\boldsymbol y}\in Y.
\label{F_cons1_safeM}
\end{equation}
Inequality~\eqref{F_cons1_safeM} can also be interpreted as a no-good interdiction cut: whenever the follower solution $\tilde{\boldsymbol y}$ attains a coverage value exceeding $t$, the leader must block at least one facility selected in $I(\tilde{\boldsymbol y})$; otherwise, the solution remains feasible and violates the target coverage requirement. Its validity is immediate: if no facility in $I(\tilde{\boldsymbol y})$ is blocked, the inequality enforces $t\geq\varphi(\tilde{\boldsymbol y})$; if at least one selected facility is blocked, the right-hand side is nonpositive and the inequality is inactive. \Fabio{This reading is exact. For a follower solution with $\varphi(\tilde{\boldsymbol y})>t$, dividing~\eqref{F_cons1_safeM} by $\varphi(\tilde{\boldsymbol y})$ gives $\sum_{i\in I(\tilde{\boldsymbol y})}x_i\geq 1-t/\varphi(\tilde{\boldsymbol y})$, and since the left-hand side is integer and the right-hand side lies in $(0,1]$, the inequality is equivalent to the no-good cut
$$
\sum_{i\in I(\tilde{\boldsymbol y})}x_i\geq 1,
$$
while for $\varphi(\tilde{\boldsymbol y})\leq t$ it is satisfied by every $\boldsymbol x\geq\boldsymbol 0$ and can be discarded. The family~\eqref{F_cons1_safeM} therefore reduces to the set-covering formulation of binary blocker problems, see Section 10.7.1 in \citet{beck2026linear} and \citet{wei2022integer}.} We use this coefficient only as a computational baseline and next exploit the max--min structure to obtain substantially stronger facility-specific coefficients.

We now derive a problem-specific facility-wise loss bound directly from the max--min structure.

\begin{proposition}
\label{prop:criticalM}
\Fabio{For each $\tilde{\boldsymbol y}\in Y$, with $Y$ defined by the upper-budget constraint~\eqref{OF_cons_l2}, taking
$M_i(\tilde{\boldsymbol y}):=\sum_{j\in J_i:\,r_j(\tilde{\boldsymbol y})=p_{ij}}d_jp_{ij}$
for every $i\in I(\tilde{\boldsymbol y})$ makes~\eqref{F_cons1} the valid family of interdiction cuts}
\begin{equation}
    t \geq 
    \varphi(\tilde{\boldsymbol y})
    -
    \sum_{i\in I(\tilde{\boldsymbol y})}
    \Fabio{\Biggl(\,\sum_{j\in J_i:\,r_j(\tilde{\boldsymbol y})=p_{ij}}d_jp_{ij}\Biggr)}
    x_i,
    \qquad \tilde{\boldsymbol y}\in Y.
\label{F_cons1_critM}
\end{equation}
\end{proposition}

The proof is given in Section~\ref{ec:proofs} of the electronic companion. In words, facility $i$ is \emph{critical} for customer $j$ when it is a least favorable selected facility reaching $j$, that is, when $r_j(\tilde{\boldsymbol y})=p_{ij}$; only such customers can lose coverage when $i$ is removed, and none of them loses more than $p_{ij}$, since $r_j$ cannot fall below zero. The loss is in fact positive only when $i$ is the sole selected facility reaching $j$: any other selected facility reaching $j$ provides at least $p_{ij}$, so removing $i$ can only raise $r_j$. For instance, if customer $j$ with $d_j=1$ is reached by two selected facilities with $p_{1j}=0.6$ and $p_{2j}=0.9$, then $r_j(\tilde{\boldsymbol y})=0.6$, facility $1$ is critical and contributes $0.6$ to the coefficient of facility $1$, facility $2$ contributes nothing, and removing facility $1$ actually raises $r_j$ to $0.9$.

\subsection{Outer exact algorithm}\label{sec:outer-algorithm}
The single-level reformulation decouples the upper-level and lower-level decisions by introducing an exponential number of inequalities associated with the set $Y$. \red{For a finite subset of these inequalities, the reformulation is an MILP.} At each iteration, a \red{relaxed master problem} contains only a subset of the interdiction cuts and produces a candidate blocker decision $\tilde{\boldsymbol x}$. The follower problem is then solved for this fixed $\tilde{\boldsymbol x}$ to determine whether there exists a follower solution whose value exceeds the target level $t$. If such a solution is found, the corresponding cut~\eqref{F_cons1} is generated for the upper-level problem. The pseudocode corresponding to this procedure is presented in Algorithm~\ref{alg1}.
\begin{algorithm}[!htbp]
\caption{Cutting-plane algorithm with interdiction cuts}
\label{alg1}
\begingroup
\setlength{\baselineskip}{0.94\baselineskip}
\algrenewcommand\algorithmicrequire{\textbf{Input:}}
\algrenewcommand\algorithmicensure{\textbf{Output:}}
\begin{algorithmic}[1]
\Require Problem data and target coverage $t$
\Ensure An optimal leader decision $\tilde{\boldsymbol x}$ and its corresponding optimal follower response $\tilde{\boldsymbol y}$
\State Initialize the \red{relaxed} upper-level problem by omitting constraints~\eqref{F_cons1}
\While{true}
    \State Solve the \red{relaxed} upper-level \red{MILP} problem and obtain an optimal blocker solution $\tilde{\boldsymbol x}$
    \State Solve the linearized lower-level problem obtained from~\eqref{OF_lower} by replacing~\eqref{OF_cons_l1} with~\eqref{linear_1}--\eqref{linear_2}, for the fixed blocker decision $\tilde{\boldsymbol x}$, and let $(\tilde{\boldsymbol y},\tilde{\boldsymbol r})$ be an optimal follower solution
    \If{\Fabio{$\xi(\tilde{\boldsymbol x}) \leq t$}}
        \State \Return $(\tilde{\boldsymbol x},\tilde{\boldsymbol y},\tilde{\boldsymbol r})$
    \EndIf
    \State Add the interdiction cut~\eqref{F_cons1} associated with $\tilde{\boldsymbol y}$ to the \red{relaxed} upper-level problem
\EndWhile
\end{algorithmic}
\endgroup
\end{algorithm}
Because $Y$ is finite and every violated follower solution generates an interdiction cut that prevents the same response from violating again, Algorithm~\ref{alg1} terminates after finitely many iterations. At termination, the current blocker solution is feasible for the complete \red{single-level reformulation} and is optimal because it is optimal for the \red{relaxed master problem}.

\section{Branch-and-Benders-Cut Algorithm for the \red{Lower-Level Problem}}\label{sec:Benders}
The follower problem is solved repeatedly within Algorithm~\ref{alg1}. Its main source of dimensionality is the customer-specific coverage vector $\boldsymbol r$ and the associated linearization constraints. This structure is particularly relevant in the \red{applications} targeted by the method, where the number of customers can be much larger than the number of candidate facilities. We therefore project out $\boldsymbol r$ and solve the follower through a branch-and-Benders-cut algorithm. The key result is that the dual separation problem admits an optimal solution in closed form, so Benders cuts can be generated directly rather than by repeatedly solving an auxiliary LP.

\red{We denote by \algname{OD+B\&BC} the outer decomposition with the lower-level problem solved by this branch-and-Benders-cut algorithm. Since the lower-level problem is itself solved by decomposition, \algname{OD+B\&BC} is a nested decomposition: it is the exact algorithm proposed in this paper.}

\subsection{Benders reformulation \red{of the lower-level problem}}\label{sec:benders-reformulation}
We exploit the mathematical structure of the linearized lower-level problem corresponding to~\eqref{OF_lower} to develop a tailored Benders decomposition algorithm. In particular, the decision variables $\boldsymbol r$ can be projected out. For a fixed blocker solution $\tilde{\boldsymbol x}$, linking constraints~\eqref{OF_cons_l3} ensure that blocked facilities cannot be selected by the follower. Therefore, we may retain the original coverage set $I_j$ in the projection of $\boldsymbol r$, since any blocked facility has $y_i=0$, does not restrict $r_j$, and contributes zero to the sum-based expression. Because the projected variables $\boldsymbol r$ appear in the objective, \emph{Benders optimality cuts} are required. Let $\vartheta\geq0$ denote the auxiliary variable representing the follower objective value. We denote by $\mathcal P$ the feasible polyhedron of the dual projection problem~\eqref{dualLP} introduced below, and by $B_\ell(\boldsymbol y,\vartheta)\geq0$ the Benders optimality cut associated with an extreme point $\ell\in\operatorname{ext}(\mathcal P)$. The resulting Benders master problem is
\begin{model}
\label{BendersOC}
\begin{align}
&& \max ~~ \vartheta \label{BendersOC_obj} \\[\modelrowsep]
&& {\rm subject~to} \quad \sum_{i\in I} f_i y_i &\leq b, & \label{BendersOC_budget} \\[\modelrowsep]
&& B_\ell(\boldsymbol y,\vartheta) &\geq 0, & \ell\in\operatorname{ext}(\mathcal P), \label{BendersOC_cuts} \\[\modelrowsep]
&& y_i &\leq 1-\tilde{x}_i, & i\in I, \label{BendersOC_linking} \\[\modelrowsep]
&& y_i &\in \{0,1\}, & i\in I, \label{BendersOC_domain_y} \\[\modelrowsep]
&& 0\leq \vartheta &\leq \bar v, & \label{BendersOC_domain_theta}
\end{align}
\end{model}
where $\bar v:=\sum_{j\in J}d_j$ is a global upper bound on the follower objective value.
\Fabio{Although the master variables are binary, the Benders cuts derived below are valid for every $\boldsymbol y\in[0,1]^m$ and therefore also strengthen the LP relaxation of the Benders master problem.}
The upper bound $\bar v$ in~\eqref{BendersOC_domain_theta}
ensures that the initial Benders master problem remains bounded even when no Benders cuts have yet been generated. Observe that any binary vector $\boldsymbol{y}$ that satisfies the budget and linking constraints of~\eqref{BendersOC} can be completed with a feasible coverage vector $\boldsymbol r$ for the linearized lower-level problem; hence, no Benders feasibility cuts are required.

For a given solution $\tilde{\boldsymbol y}\in[0,1]^m$, consider the following projection subproblem:
\begin{model}
\label{BendersSub}
\begin{align}
&& \max ~~ \sum_{j\in J} d_j r_j \label{BendersSub_obj} \\[\modelrowsep]
&& {\rm subject~to} \quad r_j &\leq 1-(1-p_{ij})\tilde y_i, & j\in J,\ i\in I_j, \label{BendersSub_cons1} \\[\modelrowsep]
&& r_j &\leq \sum_{i\in I_j} p_{ij}\tilde y_i, & j\in J, \label{BendersSub_cons2} \\[\modelrowsep]
&& r_j &\geq 0, & j\in J. \label{BendersSub_domain}
\end{align}
\end{model}
\Fabio{Problem~\eqref{BendersSub} separates by customer: each variable $r_j$ appears only in the constraints indexed by $j$, so the subproblem decomposes into $|J|$ independent single-variable problems, each maximizing $d_jr_j$ subject to two upper bounds. This separability is the structural reason why its dual admits an optimal solution in closed form.}
Let $\alpha_{ij}\geq0$ denote the dual multiplier associated with constraint~\eqref{BendersSub_cons1}, and let $\gamma_j\geq0$ denote the dual multiplier associated with constraint~\eqref{BendersSub_cons2}. The full dual of problem~\eqref{BendersSub} is
\begin{model}
\label{dualLP}
\begin{align}
&& \min ~~ \sum_{j\in J}\left[
\sum_{i\in I_j}\left(1-(1-p_{ij})\tilde y_i\right)\alpha_{ij}
+\left(\sum_{i\in I_j}p_{ij}\tilde y_i\right)\gamma_j
\right] \label{dualLP_obj} \\[\modelrowsep]
&& {\rm subject~to} \quad \sum_{i\in I_j}\alpha_{ij}+\gamma_j &\geq d_j, & j\in J, \label{dualLP_cons} \\[\modelrowsep]
&& \alpha_{ij} &\geq 0, & i\in I_j,\ j\in J, \label{dualLP_domain_alpha} \\[\modelrowsep]
&& \gamma_j &\geq 0, & j\in J. \label{dualLP_domain_gamma}
\end{align}
\end{model}

\subsection{Closed-form dual solution and Benders cuts}\label{sec:closed-form-benders}
To characterize an optimal dual solution in closed form, for a fixed solution $\tilde{\boldsymbol y}$ define the min-based and sum-based coverage terms
\begin{equation}
\label{MNj}
\tilde Q_j=
\begin{cases}
\displaystyle \min_{i\in I_j}\{1-(1-p_{ij})\tilde y_i\}, & I_j\neq\varnothing,\\
1, & I_j=\varnothing,
\end{cases}
\qquad
\tilde N_j=\sum_{i\in I_j}p_{ij}\tilde y_i,
\qquad j\in J.
\end{equation}
By~\eqref{BendersSub_cons1}--\eqref{BendersSub_cons2}, the optimal projected coverage of customer $j$ is $\min\{\tilde Q_j,\tilde N_j\}$. Define $J_Q=\{j\in J:\tilde Q_j<\tilde N_j\}$ and $J_N=\{j\in J:\tilde N_j\leq\tilde Q_j\}$. By construction, $J_Q\cap J_N=\varnothing$ and $J_Q\cup J_N=J$. For each $j\in J_Q$, select $i(j)\in\arg\min_{i\in I_j}\{1-(1-p_{ij})\tilde y_i\}$. Intuitively, $\tilde Q_j$ is the tightest of the caps~\eqref{BendersSub_cons1} imposed on $r_j$ by the selected facilities reaching $j$ (an unselected facility imposes the vacuous cap $1$), and $\tilde N_j$ is the cap~\eqref{BendersSub_cons2}, which vanishes when no selected facility reaches $j$. At a binary $\tilde{\boldsymbol y}$, $J_Q$ thus collects the customers reached by at least two selected facilities, whose coverage is dictated by the least favorable one, $i(j)$, and $J_N$ those reached by at most one. In the example of Section~\ref{sec:bigM}, $\tilde Q_j=0.6$ and $\tilde N_j=1.5$, so $j\in J_Q$ with $i(j)=1$; had only facility $1$ been selected, $\tilde Q_j=\tilde N_j=0.6$ and $j\in J_N$.
\begin{proposition}
\label{prop:closedformdual}
For a fixed solution $\tilde{\boldsymbol y}$, define the dual multipliers componentwise as
\begin{equation}
\label{closedform_alpha}
\tilde{\alpha}_{ij}
=
\begin{cases}
d_j,
& \text{if } j\in J_Q \text{ and } i=i(j),\\
0,
& \text{otherwise},
\end{cases}
\qquad j\in J,\ i\in I_j,
\end{equation}
and
\begin{equation}
\label{closedform_gamma}
\tilde{\gamma}_j
=
\begin{cases}
0,
& \text{if } j\in J_Q,\\
d_j,
& \text{if } j\in J_N,
\end{cases}
\qquad j\in J.
\end{equation}
Then $(\tilde{\boldsymbol\alpha},\tilde{\boldsymbol\gamma})$ is an optimal solution to the full dual problem~\eqref{dualLP}.
\end{proposition}

The proof is given in Section~\ref{ec:proofs} of the electronic companion.  

The feasible region of the full dual problem~\eqref{dualLP} is independent of the current solution $\tilde{\boldsymbol y}$. Hence, the multipliers defined in Proposition~\ref{prop:closedformdual} remain dual feasible when the fixed values $\tilde{\boldsymbol y}$ in the dual objective are replaced by the generic variables $\boldsymbol y$. By weak duality, this yields the following globally valid Benders optimality cut:
\begin{equation}
\underbrace{
\sum_{j\in J_Q}
d_j
\left(
1-(1-p_{i(j)j})y_{i(j)}
\right)
}_{\text{\rm min-based term}}
+
\underbrace{
\sum_{j\in J_N}
d_j
\sum_{i\in I_j}p_{ij}y_i
}_{\text{\rm sum-based term}}
\geq \vartheta.
\label{CBC}
\end{equation}
In words, $\vartheta$ is the follower value estimated by the Benders master problem, and cut~\eqref{CBC} bounds it from above by the coverage obtained when each customer in $J_Q$ is charged the cap of its critical facility $i(j)$ and each customer in $J_N$ its sum-based cap; the multipliers of Proposition~\ref{prop:closedformdual} simply place the whole weight $d_j$ of each customer on a cap that is binding at $\tilde{\boldsymbol y}$. At $\tilde{\boldsymbol y}$ the left-hand side equals $\sum_{j\in J}d_j\min\{\tilde Q_j,\tilde N_j\}$, so the cut is tight there. The quantities $\tilde Q_j$, $\tilde N_j$, and $i(j)$ can be determined by a single scan of the facilities in $I_j$ for each customer $j$. Therefore, the dual multipliers and the cut in~\eqref{CBC} can be computed in $\mathrm{O}\bigl(\sum_{j\in J}|I_j|\bigr)$ time, without solving an auxiliary LP. \Fabio{Since $\sum_{j\in J}|I_j|$ counts the facility--customer pairs with positive coverage, this is the sense in which we refer to linear-time separation.}

\subsection{Implementation}\label{sec:benders-implementation}
We implement the inner decomposition as a single-tree branch-and-Benders-cut algorithm using solver callbacks. At an incumbent solution $(\tilde{\boldsymbol y},\tilde\vartheta)$, the quantities $\tilde Q_j$ and $\tilde N_j$ are computed for all customers, the sets $J_Q$ and $J_N$ are identified, and cut~\eqref{CBC} is generated whenever its violation exceeds the tolerance $\delta$. When multiple facilities attain the minimum defining $\tilde Q_j$, we select the smallest-index minimizer; this convention is used only to make separation deterministic and reproducible. No feasibility cuts are needed because every $\boldsymbol y$ satisfying the budget and linking constraints can be completed with $\boldsymbol r=\boldsymbol 0$. The Benders cuts are valid for $\boldsymbol y\in[0,1]^m$.  We experimented with separation at fractional solutions, but the computational gains were limited; therefore, the final implementation separates cuts only at integer incumbents.

\section{Computational Study}\label{sec:DoE}
This section evaluates the two structural ingredients of the proposed method and the resulting exact algorithm. The replication code and data for this paper are available at \url{https://github.com/Yuntian-Zhang/Bilevel-MMCLBP}.

\subsection{Experimental setup and benchmark instances}\label{sec:experimental-setup}
All experiments were conducted on a machine equipped with a 3.10 GHz 8-core Intel Core i9-9900 processor and 64 GB of RAM. We used Gurobi 13.0.1 as the solver. Each run was subject to a global wall-clock time limit of 3600 seconds, including all \red{relaxed}-master and follower solves. For the \red{branch-and-Benders-cut algorithm} implemented in Gurobi, we enabled lazy constraints by setting the \texttt{LazyConstraints} parameter to 1 and set the \texttt{PreCrush} parameter to 1 to ensure that the callback-generated cuts could be properly mapped to the presolved model. Unless otherwise specified, all other Gurobi parameters were kept at their default values. The cut separation routines, including the construction of the inner Benders cuts, were implemented in C++ and integrated into the Python framework through \texttt{pybind11} to reduce the computational overhead within the callback. The C++ extension was compiled under the C++17 standard with the \texttt{-O3} optimization flag. The remaining algorithmic components were implemented in Python. The cut-separation tolerance was set to $\delta=10^{-6}$. \red{Throughout the computational study and the electronic companion, \emph{Max.\ nodes} denotes the maximum number of branch-and-bound nodes explored in any single follower solve during a complete run.}

For a fixed upper-level solution $\tilde{\boldsymbol x}$, the upper bound on each follower variable $y_i$ is set to $1-\tilde x_i$. Hence, blocked facilities are directly excluded from the follower problem. In the final configuration, the outer algorithm uses the critical-facility coefficients of Proposition~\ref{prop:criticalM}, while the follower is solved by the branch-and-Benders-cut algorithm with closed-form separation of Section~\ref{sec:Benders}.

We generate a testbed of instances following a standard spatial randomization framework commonly adopted in the covering location literature~\citep{cordeau2019benders}. Both candidate facility locations and demand points are uniformly sampled over a square region of size $30\times 30$. Demand weights are drawn independently from a discrete uniform distribution on $[1,5]$, introducing moderate heterogeneity across customers. Let $\varrho_{ij}$ denote the Euclidean distance between facility $i\in I$ and customer $j\in J$. We adopt a \red{gradual} coverage model to capture both reliable and deteriorating service effects. Let $\rho^{\mathrm{full}}$ and $\rho^{\mathrm{max}}$ denote the core and maximum covering radii, respectively, and let $\lambda>0$ denote the exponential decay parameter. Specifically, each facility provides full coverage within $\rho^{\mathrm{full}}$ and distance-dependent coverage between $\rho^{\mathrm{full}}$ and $\rho^{\mathrm{max}}$. Following~\citet{chan2016optimizing}, we set
\[
p_{ij} =
\begin{cases}
1, & \varrho_{ij} \leq \rho^{\mathrm{full}}, \\
e^{-\lambda(\varrho_{ij}-\rho^{\mathrm{full}})}, & \rho^{\mathrm{full}} < \varrho_{ij} \leq \rho^{\mathrm{max}}, \\
0, & \varrho_{ij} > \rho^{\mathrm{max}},
\end{cases}
\qquad i\in I,\ j\in J,
\]
with $\rho^{\mathrm{full}}=4$, $\rho^{\mathrm{max}}=8$, and $\lambda=\frac{1}{2}$. This setting ensures that each demand point is covered by a moderate number of facilities, while preserving sufficient variability in coverage levels. During instance generation, customers with $I_j=\varnothing$ are retained to preserve the prescribed instance size. Such customers have zero coverage under every follower solution and therefore contribute zero to both the follower objective and the generated Benders cuts.

We assume unit costs for both facility location and blocker decisions, resulting in a cardinality bound on the number of selected facilities. This setting isolates the combinatorial structure of the problem without introducing cost heterogeneity. Although more general knapsack constraints are possible, the uniform cost structure does not necessarily simplify the problem, as it increases symmetry and reduces cost-based guidance. Consequently, the computational difficulty is primarily driven by the interaction between coverage and blocker decisions. The formulation also permits nonuniform opening costs. The baseline location budget in the follower problem is set to $b=\lfloor \frac{1}{5}m\rfloor$.

To set a reasonable target coverage threshold across various settings, we first solve the unblocked follower problem using Gurobi, omitting the linking constraints~\eqref{OF_cons_l3}. The resulting problem, after linearization, is an MILP whose optimal objective value we denote by $D^0$. Let $\beta\in(0,1)$ denote the fraction of the unblocked optimum used to define the target coverage. For the baseline setting, we set $\beta=\frac{7}{10}$ and therefore $t=\beta D^0$.

\subsection{Overall performance of the final algorithm}\label{sec:overall-performance}
We first compare two variants of the same \red{outer decomposition}. Both use the critical-facility coefficients of Proposition~\ref{prop:criticalM} and differ only in how the follower problem is solved. The \red{\algname{OD+MILP}} variant, whose name records the outer decomposition and the MILP solution of the follower, solves the linearized follower problem directly with Gurobi as an MILP, whereas the \red{\algname{OD+B\&BC}} variant solves it through the inner branch-and-Benders-cut algorithm using the closed-form separation described in Section~\ref{sec:Benders}. Table~\ref{tab:benders-aggregate} summarizes the aggregate performance by instance size. For each method and each $(|I|,|J|)$ setting, we report the number of instances solved to proven optimality within the 3600-second time limit, the geometric mean solution time over the solved instances, and the average final optimality gap over all five instances, assigning a zero gap to instances solved to proven optimality. Complete instance-level results are provided in Appendix~\ref{app:detailed-results}.

\begin{table}[!t]
\centering
\small
\caption{Aggregate computational performance of the two \red{outer decomposition} variants by instance size}
\label{tab:benders-aggregate}

\setlength{\tabcolsep}{6.0pt}
\renewcommand{\arraystretch}{0.92}

\begin{tabular}{rrr rrr rrr}
\toprule
& &
& \multicolumn{3}{c}{\red{\algname{OD+MILP}}}
& \multicolumn{3}{c}{\red{\algname{OD+B\&BC}}} \\

\cmidrule(lr){4-6}
\cmidrule(lr){7-9}

$|I|$ & $|J|$ & \# inst.
& solved & \shortstack{gmean\\time (s)} & \shortstack{avg. gap\\(\%)}
& solved & \shortstack{gmean\\time (s)} & \shortstack{avg. gap\\(\%)} \\

\midrule

20 & 5000
& 5
& 5 & 126.48 & 0.00
& 5 & 4.01 & 0.00 \\

20 & 10000
& 5
& 5 & 449.87 & 0.00
& 5 & 4.29 & 0.00 \\

20 & 20000
& 5
& 5 & 1749.72 & 0.00
& 5 & 5.51 & 0.00 \\

20 & 30000
& 5
& 2 & 1617.87 & 41.00
& 5 & 6.71 & 0.00 \\

20 & 50000
& 5
& 1 & 3180.00 & 59.00
& 5 & 6.12 & 0.00 \\

\midrule

25 & 5000
& 5
& 5 & 530.60 & 0.00
& 5 & 14.93 & 0.00 \\

25 & 10000
& 5
& 4 & 1359.52 & 12.80
& 5 & 14.06 & 0.00 \\

25 & 20000
& 5
& 0 & -- & 70.40
& 5 & 15.13 & 0.00 \\

25 & 30000
& 5
& 0 & -- & 79.20
& 5 & 19.95 & 0.00 \\

25 & 50000
& 5
& 0 & -- & 88.80
& 5 & 25.64 & 0.00 \\

\midrule

35 & 5000
& 5
& 3 & 2635.31 & 27.43
& 4 & 995.46 & 14.86 \\

35 & 10000
& 5
& 0 & -- & 77.71
& 4 & 1039.85 & 15.43 \\

35 & 20000
& 5
& 0 & -- & 90.86
& 4 & 1104.86 & 15.43 \\

\bottomrule
\end{tabular}

\vspace{2pt}
\begin{minipage}{0.98\textwidth}
\footnotesize\emph{Note.}
\red{Geometric-mean times are taken row by row, over the instances that the method in that row solved; they should therefore be read together with the solved counts and the average gaps whenever the number of solved instances differs across methods.
A dash indicates that no instance in the corresponding row was solved to proven optimality.}
\end{minipage}
\end{table}

The computational advantage of the inner decomposition is substantial. \red{\algname{OD+MILP}} solves 30 of the 65 benchmark instances to proven optimality, whereas \red{\algname{OD+B\&BC}} solves 62. On the 30 instances solved by both methods, their geometric-mean solution times are 754.55 and 10.32 seconds, respectively, corresponding to a speedup of approximately 73 times. The difference becomes increasingly pronounced with problem size. For $|I|=20$, \red{\algname{OD+B\&BC}} solves all 25 instances, including all instances with 50,000 customers, whereas \red{\algname{OD+MILP}} solves 18. For $|I|=25$, \red{\algname{OD+B\&BC}} again solves all 25 instances, compared with only nine for \red{\algname{OD+MILP}}. For $|I|=35$, it proves optimality for 12 of the 15 instances, compared with only three for \red{\algname{OD+MILP}}.

The final optimality gaps reinforce this scalability advantage. Averaged over all 65 instances, the gap decreases from 42.09\% for \red{\algname{OD+MILP}} to 3.52\% for \red{\algname{OD+B\&BC}}. The contrast is particularly pronounced on the larger instances. At $(|I|,|J|)=(25,50000)$, \red{\algname{OD+MILP}} solves none of the five instances and terminates with an average gap of 88.80\%, whereas \red{\algname{OD+B\&BC}} solves all five to proven optimality. At $(|I|,|J|)=(35,20000)$, the corresponding average gaps are 90.86\% and 15.43\%, respectively. Thus, the inner decomposition improves not only the number of instances solved to proven optimality and the solution time, but also the progress achieved within the time limit on the most difficult instances.

The aggregate results in Table~\ref{tab:benders-aggregate} provide the main comparison, and Appendix~\ref{app:detailed-results} contains the complete instance-level data used to construct that summary. Overall, the results establish \red{\algname{OD+B\&BC}} as the final configuration, used in the subsequent analysis of the individual algorithmic components.

\newcommand{\DetailedResultsTable}{%
\begin{table}[!htbp]
\centering
\scriptsize
\caption{Detailed computational results for the two \red{outer decomposition} variants}
\label{tab:main-results}

\setlength{\tabcolsep}{6.pt}
\renewcommand{\arraystretch}{0.252}

\begin{tabular}{rrrrrrrrrrr}
\toprule
\multicolumn{5}{c}{Instances}
& \multicolumn{4}{c}{\red{\algname{OD+MILP}}}
& \multicolumn{2}{c}{\red{\algname{OD+B\&BC}}} \\

\cmidrule(lr){1-5}
\cmidrule(lr){6-9}
\cmidrule(lr){10-11}

$|I|$ & $|J|$ & $\rho^{\mathrm{full}}$ & $\rho^{\mathrm{max}}$ & inst.
& inc. & bound & gap (\%) & time
& obj. & time \\

\midrule

20 & 5000 & 4 & 8 & 1
& 8 & 8 & 0.00 & 196.5
& 8 & 6.8 \\

& & & & 2
& 5 & 5 & 0.00 & 31.8
& 5 & 0.9 \\

& & & & 3
& 7 & 7 & 0.00 & 192.3
& 7 & 4.4 \\

& & & & 4
& 8 & 8 & 0.00 & 235.7
& 8 & 7.1 \\

& & & & 5
& 10 & 10 & 0.00 & 114.3
& 10 & 5.4 \\

\midrule

20 & 10000 & 4 & 8 & 1
& 8 & 8 & 0.00 & 475.9
& 8 & 6.4 \\

& & & & 2
& 5 & 5 & 0.00 & 113.7
& 5 & 1.0 \\

& & & & 3
& 7 & 7 & 0.00 & 790.6
& 7 & 4.4 \\

& & & & 4
& 8 & 8 & 0.00 & 988.4
& 8 & 5.9 \\

& & & & 5
& 10 & 10 & 0.00 & 435.8
& 10 & 8.7 \\

\midrule

20 & 20000 & 4 & 8 & 1
& 9 & 9 & 0.00 & 3420.1
& 9 & 9.1 \\

& & & & 2
& 5 & 5 & 0.00 & 360.7
& 5 & 1.1 \\

& & & & 3
& 7 & 7 & 0.00 & 2628.3
& 7 & 6.8 \\

& & & & 4
& 8 & 8 & 0.00 & 3509.1
& 8 & 9.3 \\

& & & & 5
& 10 & 10 & 0.00 & 1441.4
& 10 & 8.0 \\

\midrule

20 & 30000 & 4 & 8 & 1
& 20 & 7 & 65.00 & TL
& 9 & 12.3 \\

& & & & 2
& 5 & 5 & 0.00 & 1017.3
& 5 & 1.6 \\

& & & & 3
& 20 & 5 & 75.00 & TL
& 7 & 7.3 \\

& & & & 4
& 20 & 7 & 65.00 & TL
& 8 & 8.3 \\

& & & & 5
& 10 & 10 & 0.00 & 2573.0
& 10 & 11.4 \\

\midrule

20 & 50000 & 4 & 8 & 1
& 20 & 6 & 70.00 & TL
& 8 & 8.1 \\

& & & & 2
& 5 & 5 & 0.00 & 3180.0
& 5 & 1.7 \\

& & & & 3
& 20 & 4 & 80.00 & TL
& 7 & 7.2 \\

& & & & 4
& 20 & 4 & 80.00 & TL
& 8 & 8.4 \\

& & & & 5
& 20 & 7 & 65.00 & TL
& 10 & 10.3 \\

\midrule

25 & 5000 & 4 & 8 & 1
& 9 & 9 & 0.00 & 413.8
& 9 & 12.7 \\

& & & & 2
& 7 & 7 & 0.00 & 286.1
& 7 & 6.0 \\

& & & & 3
& 9 & 9 & 0.00 & 872.7
& 9 & 23.0 \\

& & & & 4
& 8 & 8 & 0.00 & 446.1
& 8 & 12.9 \\

& & & & 5
& 11 & 11 & 0.00 & 912.5
& 11 & 32.8 \\

\midrule

25 & 10000 & 4 & 8 & 1
& 9 & 9 & 0.00 & 2081.4
& 9 & 14.4 \\

& & & & 2
& 7 & 7 & 0.00 & 740.9
& 7 & 4.7 \\

& & & & 3
& 25 & 9 & 64.00 & TL
& 9 & 30.6 \\

& & & & 4
& 9 & 9 & 0.00 & 1800.6
& 9 & 15.9 \\

& & & & 5
& 10 & 10 & 0.00 & 1230.3
& 10 & 16.7 \\

\midrule

25 & 20000 & 4 & 8 & 1
& 25 & 8 & 68.00 & TL
& 9 & 17.6 \\

& & & & 2
& 25 & 7 & 72.00 & TL
& 7 & 5.2 \\

& & & & 3
& 25 & 7 & 72.00 & TL
& 9 & 32.6 \\

& & & & 4
& 25 & 7 & 72.00 & TL
& 9 & 16.1 \\

& & & & 5
& 25 & 8 & 68.00 & TL
& 10 & 16.5 \\

\midrule

25 & 30000 & 4 & 8 & 1
& 25 & 5 & 80.00 & TL
& 9 & 22.1 \\

& & & & 2
& 25 & 6 & 76.00 & TL
& 7 & 6.1 \\

& & & & 3
& 25 & 5 & 80.00 & TL
& 9 & 38.3 \\

& & & & 4
& 25 & 4 & 84.00 & TL
& 9 & 22.5 \\

& & & & 5
& 25 & 6 & 76.00 & TL
& 11 & 27.2 \\

\midrule

25 & 50000 & 4 & 8 & 1
& 25 & 4 & 84.00 & TL
& 9 & 19.4 \\

& & & & 2
& 25 & 3 & 88.00 & TL
& 7 & 7.4 \\

& & & & 3
& 25 & 2 & 92.00 & TL
& 9 & 50.3 \\

& & & & 4
& 25 & 2 & 92.00 & TL
& 9 & 24.3 \\

& & & & 5
& 25 & 3 & 88.00 & TL
& 11 & 63.1 \\

\midrule

35 & 5000 & 4 & 8 & 1
& 35 & 13 & 62.86 & TL
& 13 & 1181.4 \\

& & & & 2
& 11 & 11 & 0.00 & 2641.5
& 11 & 1022.6 \\

& & & & 3
& 10 & 10 & 0.00 & 2492.4
& 10 & 889.5 \\

& & & & 4
& 35 & 9 & 74.29 & TL
& 35/9 & TL (74.29\%) \\

& & & & 5
& 10 & 10 & 0.00 & 2779.9
& 10 & 913.8 \\

\midrule

35 & 10000 & 4 & 8 & 1
& 35 & 8 & 77.14 & TL
& 13 & 1351.1 \\

& & & & 2
& 35 & 8 & 77.14 & TL
& 11 & 750.0 \\

& & & & 3
& 35 & 8 & 77.14 & TL
& 10 & 1167.0 \\

& & & & 4
& 35 & 6 & 82.86 & TL
& 35/8 & TL (77.14\%) \\

& & & & 5
& 35 & 9 & 74.29 & TL
& 11 & 988.7 \\

\midrule

35 & 20000 & 4 & 8 & 1
& 35 & 2 & 94.29 & TL
& 13 & 918.0 \\

& & & & 2
& 35 & 3 & 91.43 & TL
& 11 & 1005.9 \\

& & & & 3
& 35 & 4 & 88.57 & TL
& 10 & 1603.0 \\

& & & & 4
& 35 & 3 & 91.43 & TL
& 35/8 & TL (77.14\%) \\

& & & & 5
& 35 & 4 & 88.57 & TL
& 11 & 1006.7 \\

\bottomrule
\end{tabular}
\end{table}
}

\subsection{Impact of the strengthened interdiction-cut coefficients}
\label{sec:M-ablation}
We first assess the computational impact of the coefficients used in the outer interdiction cuts. The experiment is conducted on nine representative instances, consisting of three independently generated instances for each of the sizes $(|I|,|J|)=(20,20000)$, $(25,20000)$, and $(35,10000)$. These instances are used as a compact diagnostic set to isolate the contribution of the interdiction-cut coefficients while keeping all other algorithmic components fixed. \red{We compare the two choices developed in Section~\ref{sec:bigM}. The No-good variant uses the uniform solution-dependent coefficient \Fabio{$M_i(\tilde{\boldsymbol y})=\varphi(\tilde{\boldsymbol y})$, that is, the interdiction cuts~\eqref{F_cons1_safeM}}, whereas the Critical variant employs the facility-specific coefficients of Proposition~\ref{prop:criticalM}, \Fabio{that is, the interdiction cuts~\eqref{F_cons1_critM}}.} All other algorithmic components are kept unchanged, using the closed-form Benders separation. Table~\ref{tab:M-ablation} summarizes the results.

\begin{table}[!t]
\centering
\small
\caption{Effect of the interdiction-cut coefficients}
\label{tab:M-ablation}

\setlength{\tabcolsep}{6.pt}
\renewcommand{\arraystretch}{0.95}

\begin{tabular}{rrlrrrr}
\toprule
$|I|$ & $|J|$
& coefficient
& solved
& Time (s)
& \shortstack{outer\\iterations}
& \shortstack{Int.\\cuts} \\
\midrule

20 & 20000 & No-good  & 3/3 & 8.59  & 314 & 311 \\
   &       & Critical & 3/3 & 3.87  & 155 & 152 \\

\midrule

25 & 20000 & No-good  & 3/3 & 64.35 & 930 & 927 \\
   &       & Critical & 3/3 & 12.72 & 186 & 183 \\

\midrule

35 & 10000 & No-good  & 0/3 & 3600.00 & 508 & 505 \\
   &       & Critical & 3/3 & 1051.85 & 402 & 399 \\

\bottomrule
\end{tabular}

\vspace{2pt}
\begin{minipage}{0.94\textwidth}
\footnotesize\emph{Note.}
\red{``Solved'' reports the number of instances solved to proven optimality
within the time limit. Time is the geometric mean over the three runs, whereas
Outer iterations and interdiction cuts are cumulative totals. ``Int.\ cuts'' denotes the
number of interdiction cuts generated.}
\end{minipage}
\end{table}

\red{The results show a clear benefit from exploiting the facility-specific structure of the interdiction-cut coefficients. The No-good variant only enforces that at least one facility of the current follower solution be blocked, discarding any information on how much coverage each individual facility actually carries.}

\red{The Critical coefficients, in contrast, substantially strengthen the relaxed master problem. For $(|I|,|J|)=(20,20000)$, the total number of outer iterations decreases from 314 to 155 and the geometric-mean solution time from 8.59 to 3.87 seconds. For $(25,20000)$, the number of outer iterations decreases from 930 to 186 and the geometric-mean time from 64.35 to 12.72 seconds. The effect is even more pronounced for $(35,10000)$: the No-good variant solves none of the three instances within the time limit, whereas the Critical variant solves all three to proven optimality with a geometric-mean solution time of 1051.85 seconds. Since each outer iteration requires one follower solve, the reduction in iterations translates directly into fewer follower problems, which is where the running time is concentrated. These results provide direct computational evidence that the facility-specific coefficients of Proposition~\ref{prop:criticalM} are an important component of the proposed algorithm.}

\subsection{LP-based versus closed-form Benders separation}
\label{sec:lp-vs-closedform}
We next evaluate the computational benefit of exploiting the closed-form characterization of the optimal dual solution in the inner decomposition. The comparison is conducted on the same nine representative diagnostic instances introduced in Section~\ref{sec:M-ablation}. Within the same branch-and-Benders-cut algorithm, we compare \emph{LP-based separation}, which explicitly solves the Benders dual LP whenever separation is required, with \emph{closed-form separation}, which directly constructs an optimal dual solution and the corresponding cut from the analytical characterization. For a fair comparison, the LP dual model is built only once and reused across callbacks, with only its objective coefficients updated. All other algorithmic components are kept identical: both variants use the critical interdiction-cut coefficients and separate only at integer incumbent solutions. Table~\ref{tab:lp-comb-benders} reports the results.

\begin{table}[!t]
\centering
\small
\caption{LP-based versus closed-form Benders separation}
\label{tab:lp-comb-benders}
\setlength{\tabcolsep}{2.8pt}
\renewcommand{\arraystretch}{0.95}

\begin{tabular}{ccrrrrrrrr}
\toprule
$(|I|,|J|)$
& Method
& Solved
& Time
& Cuts
& Calls
& Sep. time
& ms/call
& ms/cut
& Max. nodes \\
\midrule

$(20,20000)$
& LP-based
& 3/3
& 218.22
& 8,726
& 9,768
& 878.88
& 89.98
& 100.72
& 830 \\

& Closed-form
& 3/3
& 3.64
& 6,520
& 7,527
& 4.49
& 0.60
& 0.69
& 523 \\

\midrule

$(25,20000)$
& LP-based
& 3/3
& 1167.85
& 33,795
& 35,197
& 4424.24
& 125.70
& 130.91
& 3,597 \\

& Closed-form
& 3/3
& 12.81
& 28,373
& 29,723
& 21.11
& 0.71
& 0.74
& 3,045 \\

\midrule

$(35,10000)$
& LP-based
& 0/3
& 3600.00
& 106,127
& 106,945
& 9759.54
& 91.26
& 91.96
& 87,991 \\

& Closed-form
& 3/3
& 1039.80
& 592,038
& 596,524
& 327.49
& 0.55
& 0.55
& 84,327 \\

\bottomrule
\end{tabular}

\begin{minipage}{0.99\linewidth}
\footnotesize
\emph{Note.} Time is the geometric mean over the three seeds, with time-limited runs included at their observed running time. Cuts, separation calls, and separation time are cumulative over the three seeds. The per-call and per-cut times are computed from the corresponding cumulative totals. \red{Max. nodes is a geometric mean over the three seeds.}
\end{minipage}
\end{table}

The advantage of the closed-form separation is substantial. On the six instances solved to proven optimality by both methods, the geometric-mean solution time decreases from 504.82 seconds with LP-based separation to 6.83 seconds with closed-form separation, corresponding to a speedup of approximately $73.9\times$. More importantly, the average separation time per callback decreases from 117.94 ms to 0.687 ms, a reduction by a factor of approximately $171.6$, while the corresponding per-cut time decreases from 124.72 ms to 0.734 ms. Despite reusing the dual model, repeated dual-LP reoptimization still accounts for a substantial share of the cost of LP-based separation.

The difference also translates into substantially better scalability. Closed-form separation solves all nine instances, whereas LP-based separation solves only six. In particular, the LP-based variant reaches the 3600-second time limit on all three $(35,10000)$ instances, with final optimality gaps of 74.3\%, 77.1\%, and 77.1\%, while the closed-form variant proves optimality for all three. At the same time, the maximum branch-and-bound tree sizes are broadly comparable, indicating that the improvement does not primarily arise from systematically smaller follower search trees. Rather, exploiting the problem structure eliminates repeated LP reoptimization and reduces the cost of Benders separation by roughly two orders of magnitude, which translates directly into substantial end-to-end gains.

\subsection{Scalability with respect to the customer dimension}
\label{sec:scalability}
The inner decomposition is designed for instances in which the number of candidate facilities is moderate while the customer dimension can be substantially larger. The main benchmark already suggests this behavior: for $|I|=20$ and $25$, \red{\algname{OD+B\&BC}} solves every tested instance with up to 50,000 customers, whereas \red{\algname{OD+MILP}} deteriorates rapidly as $|J|$ increases. To examine this regime more systematically, we conduct an additional scalability experiment with $|I|\in\{20,25,35\}$ and $|J|\in\{100000,150000,200000,250000,300000\}$. For each combination, we generate three independent instances using the same data-generation procedure as in the main benchmark. We fix $\beta=0.7$ and $b=\lfloor0.2|I|\rfloor$ and use the final \red{\algname{OD+B\&BC}} configuration. Each run is subject to the same time limit of 3600 seconds.

Table~\ref{tab:scalability-customer} reports the results. The algorithm exhibits strong scalability with respect to the customer dimension when the number of candidate facilities remains moderate. For both $|I|=20$ and $|I|=25$, all 15 instances are solved to proven optimality, including all instances with 300,000 customers. For $|I|=20$, the geometric-mean solution time increases from 7.09 seconds at 100,000 customers to only 16.31 seconds at 300,000 customers. For $|I|=25$, the corresponding increase is from 34.93 to 80.63 seconds. Moreover, neither the number of outer iterations nor the maximum follower branch-and-bound tree exhibits a systematic increase with $|J|$ for these two facility sizes. This indicates that increasing the customer dimension mainly increases the cost of processing a follower problem, rather than fundamentally changing the combinatorial difficulty of \red{the outer decomposition}.

\begin{table}[htbp]
\centering
\small
\caption{Scalability with respect to the customer dimension}
\label{tab:scalability-customer}
\setlength{\tabcolsep}{6.pt}
\renewcommand{\arraystretch}{0.95}
\begin{tabular}{rrrrrrr}
\toprule
$|I|$ & $|J|$ & Solved & Time (s) & Gap (\%) &
Outer iter. & Max. nodes \\
\midrule

\multirow{5}{*}{20}
& 100,000 & 3/3 & 7.09  & -- & 35.9 & 518 \\
& 150,000 & 3/3 & 8.98  & -- & 33.0 & 670 \\
& 200,000 & 3/3 & 13.13 & -- & 39.7 & 637 \\
& 250,000 & 3/3 & 15.62 & -- & 39.6 & 520 \\
& 300,000 & 3/3 & 16.31 & -- & 36.7 & 509 \\

\midrule

\multirow{5}{*}{25}
& 100,000 & 3/3 & 34.93 & -- & 63.9 & 3,485 \\
& 150,000 & 3/3 & 50.04 & -- & 64.3 & 3,124 \\
& 200,000 & 3/3 & 57.22 & -- & 60.0 & 2,974 \\
& 250,000 & 3/3 & 76.89 & -- & 68.3 & 3,468 \\
& 300,000 & 3/3 & 80.63 & -- & 57.3 & 3,339 \\

\midrule

\multirow{5}{*}{35}
& 100,000 & 3/3 & 2174.00 & --    & 141.3 & 80,500 \\
& 150,000 & 3/3 & 2445.15 & --    & 127.9 & 80,627 \\
& 200,000 & 2/3 & 2272.57 & 62.86 & 111.6 & 85,923 \\
& 250,000 & 2/3 & 2538.66 & 62.86 & 113.7 & 76,844 \\
& 300,000 & 1/3 & 3180.50 & 70.00 & 112.1 & 79,093 \\

\bottomrule
\end{tabular}

\vspace{2pt}
\begin{minipage}{0.98\linewidth}
\footnotesize
\emph{Note.}
Time is the geometric mean over instances solved to proven optimality.
Gap is the arithmetic mean final optimality gap over time-limited instances
and is omitted when all three instances are solved.
Outer iter. and Max. nodes are geometric means over all three runs.
\end{minipage}
\end{table}

The picture changes markedly when the facility dimension increases to $|I|=35$. All instances with 100,000 and 150,000 customers are still solved to optimality, but the geometric-mean running times already exceed 2,000 seconds. At 200,000 and 250,000 customers, two of the three instances are solved within the time limit, while only one instance is solved at 300,000 customers. The average final gaps of the unsolved runs are 62.86\%, 62.86\%, and 70.00\%, respectively. At the same time, the maximum follower branch-and-bound tree increases from only a few thousand nodes for $|I|=25$ to approximately $8\times10^4$ nodes for $|I|=35$, even though the number of outer iterations remains of a comparable order across the tested customer sizes. These results indicate that the principal computational bottleneck eventually shifts to the combinatorial complexity of the follower problem as the number of candidate facilities increases.

Overall, the dedicated experiment extends the favorable customer-scaling regime observed in the main benchmark from 50,000 to 300,000 customers. The results therefore support scalability primarily with respect to the customer dimension when the number of candidate facilities is moderate, rather than a generic claim of large-scale scalability across both problem dimensions.

\subsection{Virginia Beach AED  case study}\label{sec:case-study}
To illustrate the practical relevance of the MMCLBP, we use the Virginia Beach dataset introduced by~\citet{custodio2022spatiotemporal}, available at \url{https://github.com/INFORMSJoC/2020.1022}. The dataset contains 2,706 geocoded out-of-hospital cardiac-arrest (OHCA) events observed from January 1, 2017 to June 30, 2019 and 40 emergency medical service (EMS), fire, and police station locations. We use the OHCA events as customers and the emergency-service facilities as candidate AED  bases. In an operational AED  system, however, the base that ultimately serves an incident may depend on real-time device availability and single-responder capabilities, 
and is therefore not fully controlled by the planner. We accordingly interpret the selected bases as an admissible service pool and use the minimum coverage level over that pool as a conservative guarantee against an unfavorable realized service origin. Each event has unit weight, $d_j=1$.

\orange{The model does not optimize the dispatch decision: it evaluates service through the least favorable admissible origin. A base that lowers this guarantee therefore does not necessarily worsen the performance of an operational dispatch system, and the criterion does not represent a dispatcher that can always select the most favorable available base.}

For each base--OHCA pair, the coverage level is computed from Haversine distance using the same functional form as in the synthetic experiments. We set $\rho^{\mathrm{full}}=4$ km, $\rho^{\mathrm{max}}=8$ km, and $\lambda=1/2$. The 4-km radius approximately matches the 97th percentile of the distance from an OHCA event to its nearest candidate base, while 2,705 of the 2,706 events have at least one base within 8 km. In the absence of reliable heterogeneous cost data, we set $f_i=e_i=1$ for every candidate base. We use $b=10$ and $\beta=0.7$, with target $t=\beta D^0$ defined from the optimal unblocked follower value.

\begin{figure}[!htbp]
\centering
\includegraphics[width=0.48\textwidth]{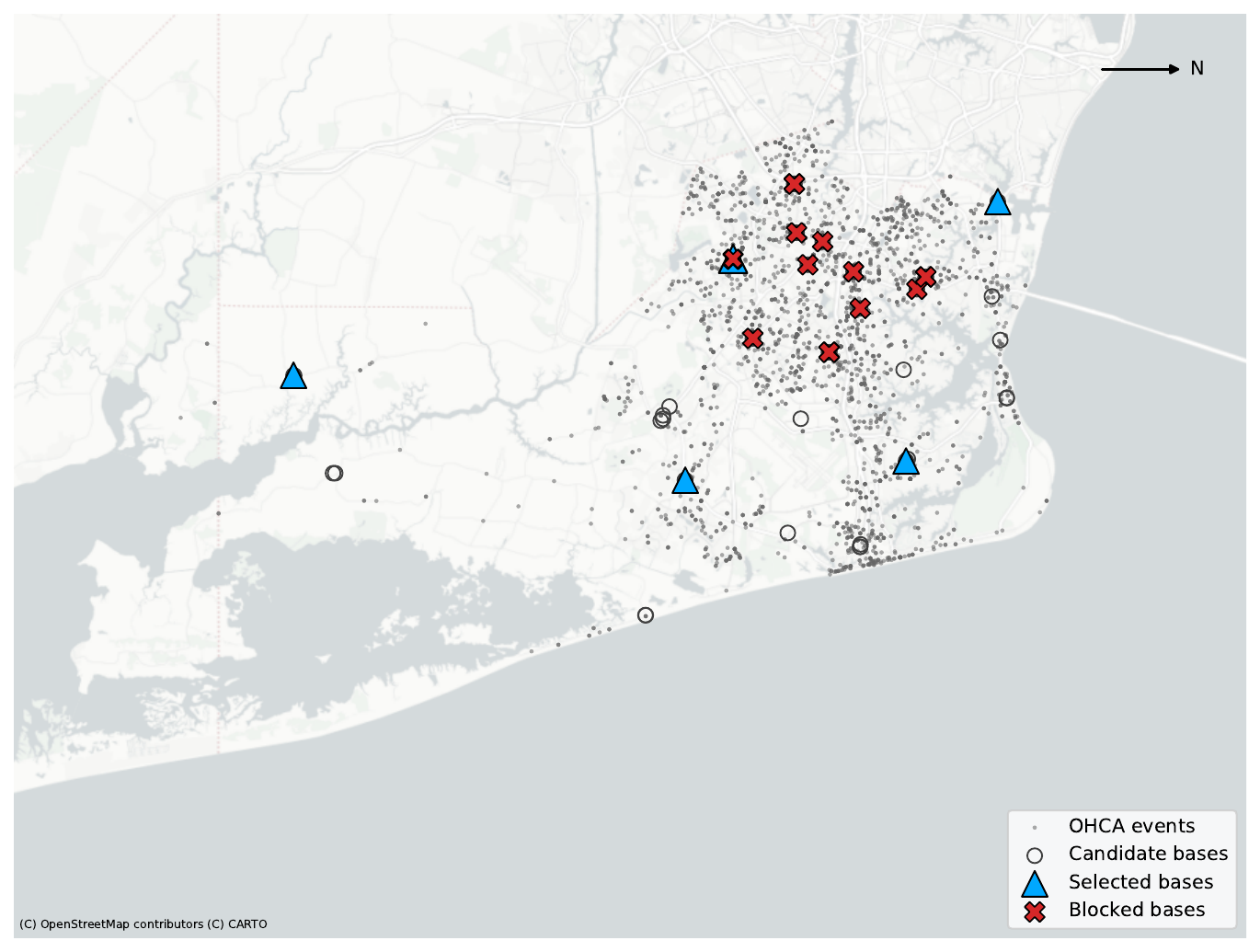}
\caption{Spatial distribution of OHCA events, candidate bases, selected bases, and optimal blocking locations in Virginia Beach. Four pairs coincide geographically, so 15 blocked bases appear as 11 markers.}
\label{fig:VB}
\end{figure}

Under this max--min criterion, the unblocked follower activates only five bases, $\{1,7,15,19,38\}$, even though the budget allows up to ten. \orange{Because the objective evaluates the least favorable admissible service origin, deploying an additional base need not improve the value of a fixed configuration. This explains why an optimal deployment may leave part of the budget unused. Nevertheless, increasing the budget cannot decrease the optimal follower value, since the budget constraint is an upper bound and the previous deployment remains feasible.} The unblocked follower value is $D^0=1775.42$ and the target is $t=\beta D^0=1242.80$. The optimal blocker objective is 15. Because each blockage has unit cost, this has a direct resilience interpretation: the loss of any at most 14 candidate bases leaves the optimal guaranteed coverage strictly above the target, whereas a simultaneous loss of the 15 reported blocking bases is sufficient to reduce it to the target or below. The optimal blocking set is $\{4,8,10,11,12,18,20,21,24,26,31,32,35,36,38\}$. Figure~\ref{fig:VB} shows the two sets against the spatial distribution of the OHCA events. Strikingly, the nominal deployment and the blocker set overlap in the single location 38. Thus, locations that are attractive in nominal operation need not be those that are most important for preserving the follower's reconfiguration flexibility after disruptions.

The blocking locations must be interpreted jointly: they form a minimum disruption set whose simultaneous removal is sufficient to restrict the follower below the prescribed target, \red{no matter where the remaining up to $b$ facilities are open.  The blocking locations } are not claimed to be individually critical when removed in isolation. \orange{This is a statement about the modeled coverage guarantee, rather than a claim that the same disruptions necessarily drive the performance of an operational dispatch policy below the target.} From a resilience perspective, the result illustrates why protecting only the facilities selected in a nominal deployment can be misleading: locations that are unused nominally may still provide valuable reconfiguration options after other facilities are disrupted.

From a managerial perspective, location 38 is a natural priority for protection or reinforcement: it is the only base that is both selected in the nominal deployment and included in the minimum blocking set. If it is unavailable together with the other 14 locations in that set, the decision maker cannot guarantee the target coverage, regardless of how the surviving bases are configured. This does not imply that protecting location 38 alone is sufficient; rather, it identifies the location as a priority for redundancy and continuity measures within a joint resilience strategy.

\section{Conclusions}\label{sec:conclusions}
We introduced the Max--Min Covering Location Blocker Problem, a bilevel resilience model for covering systems in which guaranteed customer service is determined by the least favorable selected facility capable of serving each customer. The solution approach exploits this structure at both levels of a nested exact decomposition. Facility-specific loss coefficients strengthen the outer interdiction-cut approximation, while a closed-form optimal dual solution enables Benders optimality cuts for the follower to be generated in closed form without auxiliary LP solves. The computational experiments show that these two ingredients are responsible for the main performance gains: the strengthened outer coefficients sharply reduce repeated follower solves, and closed-form separation reduces the cost of Benders cut generation by roughly two orders of magnitude. On the current benchmark, \red{\algname{OD+B\&BC}} solves 62 of 65 instances to proven optimality, compared with 30 for \red{\algname{OD+MILP}}. Against two general-purpose bilevel solvers from the literature, run in the configurations applicable to the MMCLBP, the proposed method is one to three orders of magnitude faster; it proves optimality on every instance of the comparison, including those that neither solver closes within the time limit, and needs a few dozen follower solves where they explore tens of thousands of nodes. Appendix~\ref{app:general-purpose} reports this comparison.

The Virginia Beach case study \orange{illustrates the blocker model under the max--min service criterion adopted here, and} shows that nominally selected facilities and disruption-critical locations can be very different, highlighting the role of reconfiguration flexibility in resilience analysis. Additional scalability experiments show that the final \red{\algname{OD+B\&BC}} algorithm solves all tested instances with $|I|\in\{20,25\}$ and up to 300,000 customers to proven optimality, while the harder $|I|=35$ instances indicate that the combinatorial complexity of the follower eventually becomes the dominant bottleneck as the facility dimension grows.

\bibliographystyle{informs2014}
\bibliography{biblio} 

\newpage
\begin{APPENDICES}
\section{Comparison with general-purpose bilevel solvers}
\phantomsection
\label{app:general-purpose}
The MMCLBP admits a formulation as a single mixed-integer bilevel linear
program and can therefore be submitted to a general-purpose bilevel solver, or
addressed through generic network and value-function relaxations such as those
of \citet{lozano2026network}.
The purpose of this appendix is to establish how such solvers perform on the
problem, and thereby to assess whether the tailored decomposition of
Sections~\ref{sec:outer-decomposition} and~\ref{sec:Benders} is warranted. We
describe the two solvers considered and the single-level program they receive,
specify the test instances, report the computational results, and close with
an interpretation of the outcome.

We consider two general-purpose solvers for mixed-integer bilevel linear
programs: the intersection-cut branch-and-cut algorithm
of~\citet{fischetti2017general,fischetti2018intersection}, in the binary
distributed by its authors, and version 1.2.2 of the open-source solver MibS
of~\citet{tahernejad2020branch}. Both read the same input format, in which
the program is supplied as its high-point relaxation together with a file
declaring which variables and constraints belong to the lower level and what
the lower-level objective is. Both are run single-threaded and in
deterministic mode.

Both solvers receive the bilevel program obtained by replacing
constraint~\eqref{OF_cons_l1} of the lower-level problem~\eqref{OF_lower} by
the linear inequalities~\eqref{linear_1}--\eqref{linear_2}, whose constraints
are then linear at both levels. Its high-point relaxation is
\begin{equation}
\min\Big\{\sum_{i\in I}e_ix_i \;:\;
\sum_{j\in J}d_jr_j\leq t,\;
\eqref{linear_1},\;\eqref{linear_2},\;\eqref{OF_cons_l2},\;
y_i+x_i\leq1,\; i\in I,\;
\boldsymbol x,\boldsymbol y\in\{0,1\}^m,\;
\boldsymbol r\in[0,1]^n\Big\},
\label{eq:hpr}
\end{equation}
the follower maximizing $\sum_{j\in J}d_jr_j$ over the same constraints with
the target-coverage requirement $\sum_{j\in J}d_jr_j\leq t$ omitted.  

The first solver offers several families of bilevel cuts, and we use the
hypercube family, denoted HC++, which is the family applicable to the MMCLBP.
The remaining families require the lower-level problem to be purely integral,
in variables as well as in constraint coefficients, whereas the coverage
variables $\boldsymbol r$ are continuous and the coefficients $p_{ij}$ and
$1-p_{ij}$ in~\eqref{linear_1}--\eqref{linear_2} are fractional.  

The instances are generated by the procedure of
Section~\ref{sec:experimental-setup} on the factorial grid
$|I|\in\{10,15,20\}$ by $|J|\in\{50,500,5{,}000\}$, one instance per cell.
The three values of $|I|$ give facility budgets $b=2,3,4$, and the largest
cell, $|I|=20$ and $|J|=5{,}000$, is the smallest size in the benchmark of
Section~\ref{sec:DoE}; beyond this range a general-purpose solver cannot be
applied at all. For $|I|=10$ the optimal value is obtained by complete
enumeration of all $2^{|I|}$ blocker decisions, which also validates the
single-level program~\eqref{eq:hpr}; for $|I|=15$ and $|I|=20$ it is obtained
by the proposed algorithm.

Table~\ref{tab:gp-comparison} reports, for each cell of the grid, the optimal
value, the computing time of the three methods, the slowdown of each
general-purpose solver relative to the \red{\algname{OD+B\&BC}} algorithm, and
the search effort each method expended. Search effort is measured by the
quantity natural to each method: for \red{\algname{OD+B\&BC}} it is the number of outer
iterations of Algorithm~\ref{alg1}; each iteration solves the follower problem
once, and each nonterminal iteration generates one interdiction cut. For HC++ and MibS, both of
which are branch-and-cut algorithms, it is the number of branch-and-bound nodes processed.
Search effort is reported only where the method proved optimality, since at
the time limit it measures the work performed within 600 seconds rather than
the work required. All runs are single-threaded, with a time limit of 600
seconds, on the machine of Section~\ref{sec:experimental-setup}; the times of
the proposed algorithm are medians of three runs on the same machine, so that
the three time columns are directly comparable.

\begin{table}[!htbp]
\centering
\caption{General-purpose bilevel solvers and the proposed algorithm}
\label{tab:gp-comparison}

\setlength{\tabcolsep}{5.0pt}
\renewcommand{\arraystretch}{0.92}
\small

\begin{tabular}{rrr rrr rr rrr}
\toprule
& & & \multicolumn{3}{c}{Time (s)} & \multicolumn{2}{c}{Slowdown factor} & \multicolumn{3}{c}{Search effort} \\
\cmidrule(lr){4-6}\cmidrule(lr){7-8}\cmidrule(lr){9-11}
$|I|$ & $|J|$ & obj. & \red{\algname{OD+B\&BC}} & HC++ & MibS & HC++ & MibS & \red{\algname{OD+B\&BC}} & HC++ & MibS \\
\midrule
10 & 50 & 4 & 0.02 & 0.27 & 1.25 & 16 & 74 & 12 & 153 & 1,393 \\
10 & 500 & 6 & 0.04 & 3.06 & 55.91 & 81 & 1,471 & 18 & 402 & 10,669 \\
10 & 5,000 & 7 & 0.07 & 23.20 & TL & 341 & -- & 24 & 140 & -- \\
15 & 50 & 6 & 0.11 & 3.47 & 120.25 & 32 & 1,124 & 34 & 1,462 & 125,435 \\
15 & 500 & 7 & 0.36 & 60.10 & TL & 165 & -- & 38 & 2,513 & -- \\
15 & 5,000 & 8 & 0.90 & TL & TL & -- & -- & 67 & -- & -- \\
20 & 50 & 5 & 0.78 & 394.98 & TL & 510 & -- & 32 & 70,627 & -- \\
20 & 500 & 6 & 2.06 & TL & TL & -- & -- & 51 & -- & -- \\
20 & 5,000 & 8 & 3.26 & TL & TL & -- & -- & 66 & -- & -- \\
\bottomrule
\end{tabular}

\vspace{2pt}
\begin{minipage}{0.98\linewidth}
\footnotesize
\emph{Note.}
``obj.'' is the optimal number of blocked locations. \red{TL indicates a run
terminated at the time limit. Slowdown factors are relative to \red{\algname{OD+B\&BC}} on the
same instance and are reported only for instances solved to proven optimality
by the general-purpose solver.}
\end{minipage}
\end{table}

Both general-purpose solvers return the optimal value whenever they prove
optimality. HC++ solves six of the nine cells and MibS three, whereas the
proposed algorithm solves all nine in at most $3.26$ seconds. On the cells
they close, the slowdown grows monotonically in both dimensions, from $16$ to
$510$ for HC++ and from $74$ to $1{,}471$ for MibS. The failure of MibS at
$(|I|,|J|)=(20,50)$ shows that the facility dimension alone can defeat a
general-purpose solver even with only fifty customers.

Two observations explain this behavior.

First, relaxation~\eqref{eq:hpr} provides no information. Its optimal value is
$0$, attained at $\boldsymbol x=\boldsymbol 0$, $\boldsymbol y=\boldsymbol 0$,
$\boldsymbol r=\boldsymbol 0$, since the leader may block no location while
the follower is not yet required to respond. The root bound available to a
general-purpose solver is therefore $0$ on every instance, and the search
relies entirely on bilevel-feasibility cuts. After its own root processing,
HC++ attains a root bound of $2$ for $|I|=10$ and of $1$ for $|I|\geq15$,
against optimal values between $4$ and $8$. The consequence is visible in the
last three columns of Table~\ref{tab:gp-comparison}: closing $|I|=20$,
$|J|=50$ costs HC++ $70{,}627$ nodes and defeats MibS, whereas the proposed
algorithm needs $32$ outer iterations, since each nonterminal iteration adds an
interdiction cut whose coefficients come from
Proposition~\ref{prop:criticalM} rather than a branching decision.

Second, a general-purpose solver must be supplied with
program~\eqref{eq:hpr} in its entirety, so that the size of the program grows
with the number of customers. At $|I|=20$ and $|J|=5{,}000$ the lower-level
problem already contributes $5{,}020$ variables and $21{,}992$ constraints,
and at $|J|=300{,}000$ it contributes more than $1.3$ million constraints. The
proposed algorithm never constructs this program: the upper-level problem is
maintained as a \red{relaxed master} over $m$ binary variables, and the coverage
variables $\boldsymbol r$ are projected out of the follower by the Benders
reformulation of Section~\ref{sec:Benders}.


\end{APPENDICES}

\ECSwitch
\ECHead{Electronic Companion}

\section{Notation and data}
\phantomsection
\label{ec:notation}
\Fabio{This section collects the notation and data definitions used throughout the manuscript. Sets, indices, input data, decisions, and decomposition symbols are organized by category in Table~\ref{tab:notation}. Domains and set memberships are stated next to each symbol to keep the model formulations concise and avoid repeating these definitions in the main text.}
\begingroup
\renewcommand{\Fabio}[1]{#1}
\small
\renewcommand{\arraystretch}{1.06}
\setlength{\LTleft}{0pt}
\setlength{\LTright}{\fill}
\begin{longtable}{@{}p{0.30\textwidth}p{0.66\textwidth}@{}}
\caption{Notation and data}\label{tab:notation}\\
\toprule
\textbf{Notation} & \textbf{Description} \\
\midrule
\endfirsthead
\multicolumn{2}{@{}l}{\small\itshape Table~\thetable\ (continued)}\\
\toprule
\textbf{Notation} & \textbf{Description} \\
\midrule
\endhead
\midrule
\multicolumn{2}{r@{}}{\small\itshape Continued on next page}\\
\endfoot
\bottomrule
\endlastfoot
\multicolumn{2}{@{}l}{\textbf{Sets and indices}} \\
\addlinespace
$m\in\mathbb{N}$ & Number of candidate facility locations. \\
$I=\{1,2,\dots,m\}$ & Set of candidate facility locations. \\
$n\in\mathbb{N}$ & Number of customers. \\
$J=\{1,2,\dots,n\}$ & Set of customers. \\
\midrule
\multicolumn{2}{@{}l}{\textbf{Input data and parameters}} \\
\addlinespace
\Fabio{$f_i\in\mathbb{R}_{\geq 0}$} & \Fabio{Facility-opening cost at location $i$.} \\
\Fabio{$e_i\in\mathbb{R}_{\geq 0}$} & \Fabio{Cost of blocking location $i$.} \\
\Fabio{$d_j\in\mathbb{R}_{> 0}$} & \Fabio{Importance weight of customer $j$.} \\
\Fabio{$p_{ij}\in[0,1]$} & \Fabio{Coverage level provided by location $i$ to customer $j$.} \\
$I_j=\{i\in I:0<p_{ij}\leq 1\}$ & Locations providing positive coverage to customer $j$. \\
$J_i=\{j\in J:0<p_{ij}\leq 1\}$ & Customers receiving positive coverage from location $i$. \\
\Fabio{$b\in\mathbb{R}_{\geq 0}$} & \Fabio{Facility-opening budget.} \\
\Fabio{$t\in\mathbb{R}_{\geq 0}$} & \Fabio{Target coverage threshold.} \\
\midrule
\multicolumn{2}{@{}l}{\textbf{Decisions and coverage}} \\
\addlinespace
$\boldsymbol{x}\in\{0,1\}^m$ & Blocker decision vector; $x_i=1$ if location $i$ is blocked and $x_i=0$ otherwise. \\
$\boldsymbol{y}\in\{0,1\}^m$ & Facility-location decision vector; $y_i=1$ if location $i$ is selected and $y_i=0$ otherwise. \\
\Fabio{$r_j\in\mathbb{R}_{\geq 0}$} & \Fabio{Overall coverage level of customer $j$.} \\
$I(\boldsymbol{y})=\{i\in I:y_i=1\}$ & Locations selected under $\boldsymbol{y}$. \\
$I_j(\boldsymbol{y})=I_j\cap I(\boldsymbol{y})$ & Selected locations providing positive coverage to customer $j$. \\
$Y=\{\boldsymbol{y}\in\{0,1\}^m:\sum_{i\in I}f_iy_i\leq b\}$ & Follower-feasible facility-opening set before blocking. \\
$\varphi(\boldsymbol{y})=\sum_{j\in J}d_jr_j(\boldsymbol{y})$ & Total worst-case partial coverage under $\boldsymbol{y}$. \\
$\xi(\boldsymbol{x})\in\mathbb{R}_{\geq 0}$ & \red{Optimal-value function of the lower-level problem: optimal follower value under blocker decision $\boldsymbol{x}$.} \\
\midrule
\multicolumn{2}{@{}l}{\textbf{Decomposition notation}} \\
\addlinespace
$\tilde Q_j\in[0,1]$ & Min-based coverage term at the fixed solution $\tilde{\boldsymbol y}$; set to $1$ if $I_j=\varnothing$. \\
$\tilde N_j\in\mathbb{R}_{\geq 0}$ & Sum-based coverage term at the fixed solution $\tilde{\boldsymbol y}$; the empty sum equals $0$. \\
\Fabio{$\bar v=\sum_{j\in J}d_j$} & \Fabio{Global upper bound on the follower objective value.} \\
\Fabio{$\vartheta\in[0,\bar v]$} & \Fabio{Follower-value variable in the Benders master problem.} \\
\Fabio{$\mathcal P$} & \Fabio{Feasible polyhedron of the full dual Benders subproblem.} \\
$J_Q=\{j\in J:\tilde Q_j<\tilde N_j\}$ & Customers using the min-based term at $\tilde{\boldsymbol y}$. \\
$J_N=\{j\in J:\tilde N_j\leq\tilde Q_j\}$ & Customers using the sum-based term at $\tilde{\boldsymbol y}$. \\
$i(j)$ & Minimizing facility index used in the min-based Benders term for customer $j$. \\
$\alpha_{ij},\gamma_j\geq0$ & Dual multipliers associated with~\eqref{BendersSub_cons1} and~\eqref{BendersSub_cons2}, respectively. \\
\midrule
\multicolumn{2}{@{}l}{\textbf{Interdiction-cut loss coefficients}} \\
\addlinespace
$M_i(\tilde{\boldsymbol y})\geq0$ & Generic valid loss coefficient for facility $i$ in the interdiction cut. \\
\Fabio{$M_i(\tilde{\boldsymbol y})=\varphi(\tilde{\boldsymbol y})$} & \Fabio{Uniform coefficient, yielding the no-good interdiction cuts~\eqref{F_cons1_safeM}.} \\
\Fabio{$M_i(\tilde{\boldsymbol y})=\sum_{j\in J_i:\,r_j(\tilde{\boldsymbol y})=p_{ij}}d_jp_{ij}$} & \Fabio{Critical-facility coefficient, yielding the interdiction cuts~\eqref{F_cons1_critM}.} \\
$M_i^-(\tilde{\boldsymbol y})$ & Shorthand used in the proof of Proposition~\ref{prop:criticalM} for the critical-facility coefficient. \\
$S\subseteq I(\tilde{\boldsymbol y})$ & Subset of selected facilities blocked simultaneously. \\
$\tilde{\boldsymbol y}^{\,S}$ & Surviving follower solution obtained from $\tilde{\boldsymbol y}$ by deleting the facilities in $S$. \\
\midrule
\multicolumn{2}{@{}l}{\textbf{Computational notation}} \\
\addlinespace
$\delta>0$ & Cut-separation tolerance in the callback. \\
$\varrho_{ij}$ & Euclidean distance between facility $i$ and customer $j$. \\
$\rho^{\mathrm{full}},\rho^{\mathrm{max}}$ & Core and maximum coverage radii used in instance generation. \\
$\lambda>0$ & Exponential coverage-decay parameter. \\
$\beta\in(0,1)$ & Fraction of the unblocked optimum used to define the target $t$. \\
$D^0$ & Optimal follower value of the unblocked problem, used to set $t=\beta D^0$. \\
\end{longtable}
\endgroup

\section{\red{Non-submodularity of the coverage functions}}
\phantomsection
\label{ec:submodularity}
\red{Recall that a set function $g$ on a ground set $E$ is submodular when
$g(A\cup\{e\})-g(A)\geq g(B\cup\{e\})-g(B)$ for all $A\subseteq B\subseteq E$ and
$e\in E\setminus B$. Two functions are relevant here, and they are defined on different
ground sets: the follower objective $\varphi$, as a function of the set of selected
locations, and the optimal follower value $\xi$, as a function of the set of blocked
locations. Neither is submodular, as the following two instances show. Both are stated
with unit facility-opening and blocking costs.}

\begin{example}
\label{ex:phi-nonsubmodular}
\red{One customer with $d_1=1$ and three candidate locations with
$p_{11}=0.4$, $p_{21}=0.2$ and $p_{31}=0.3$. Take $A=\{1\}$, $B=\{1,2\}$ and $e=3$. Then
$\varphi(A)=0.4$ and $\varphi(A\cup\{e\})=\min\{0.4,0.3\}=0.3$, whereas
$\varphi(B)=\min\{0.4,0.2\}=0.2$ and $\varphi(B\cup\{e\})=\min\{0.4,0.2,0.3\}=0.2$. The
marginal contribution of location $3$ is $-0.1$ on $A$ and $0$ on $B$, so submodularity
fails. The same values show that $\varphi$ is not monotone, since
$\varphi(\{1\})>\varphi(\{1,2\})$.}
\end{example}

\begin{example}
\label{ex:xi-nonsubmodular}
\red{Two customers with $d_1=d_2=1$, three candidate locations, and follower budget $b=2$.
The coverage levels are $(p_{11},p_{12})=(0.2,0)$, $(p_{21},p_{22})=(0.4,0.4)$ and
$(p_{31},p_{32})=(0,1)$, so that location $1$ can serve only customer $1$ and location $3$
only customer $2$. Unblocked, the follower opens $\{1,3\}$ and attains
$\xi(\varnothing)=0.2+1=1.2$; opening $\{2,3\}$ instead would give only $0.8$, because
location $2$ enters the minimum of customer $2$ and lowers its guaranteed coverage from
$1$ to $0.4$. Blocking gives
$$
\xi(\{3\})=0.8,\qquad \xi(\{1\})=1,\qquad \xi(\{1,3\})=0.8,
$$
attained by $\{2\}$, $\{3\}$ and $\{2\}$, respectively. Take $A=\varnothing$, $B=\{1\}$
and $e=3$: the marginal contribution of blocking location $3$ is $-0.4$ on $A$ and $-0.2$
on $B$, and submodularity would require $-0.4\geq-0.2$.}
\end{example}

\red{The last example  shows that submodularity cannot be assumed for the MMCLBP, which is
what the argument of Section~\ref{sec:literature} requires: the framework
of~\citet{taninmis2022branch} applies to followers maximizing a monotone submodular set
function, and is therefore unavailable here. This is why the loss coefficients of
Proposition~\ref{prop:criticalM} are derived from the max--min structure rather than from
marginal gains.}

\section{Proofs}
\phantomsection
\label{ec:proofs}
This section contains the proofs omitted from Sections~\ref{sec:vf-reformulation},~\ref{sec:bigM}, and~\ref{sec:closed-form-benders}. The statements are repeated for convenience.

\begin{repeatproposition}
\red{Let the coefficients $M_i(\tilde{\boldsymbol y})$ satisfy
$$
\varphi(\tilde{\boldsymbol y}^{\,S})
\geq
\varphi(\tilde{\boldsymbol y})-\sum_{i\in S}M_i(\tilde{\boldsymbol y}),
\qquad
\tilde{\boldsymbol y}\in Y,\ S\subseteq I(\tilde{\boldsymbol y}),
$$
where $\tilde{\boldsymbol y}^{\,S}$ is obtained from $\tilde{\boldsymbol y}$ by setting
$\tilde y_i=0$ for $i\in S$. Then $\boldsymbol x\in\{0,1\}^m$ satisfies every
inequality~\eqref{F_cons1} if and only if $t\geq\xi(\boldsymbol x)$.}
\end{repeatproposition}

\begin{myproof}
\red{Assume first that $t\geq\xi(\boldsymbol x)$ and fix $\tilde{\boldsymbol y}\in Y$. Let
$S(\boldsymbol x)=\{i\in I(\tilde{\boldsymbol y}):x_i=1\}$. The surviving solution
$\tilde{\boldsymbol y}^{\,S(\boldsymbol x)}$ is feasible for the follower under
$\boldsymbol x$: it satisfies the budget~\eqref{OF_cons_l2} because
$\tilde{\boldsymbol y}^{\,S(\boldsymbol x)}\leq\tilde{\boldsymbol y}$ componentwise and
$f_i\geq0$, and it satisfies the linking constraints~\eqref{OF_cons_l3} because
$\tilde y_i^{\,S(\boldsymbol x)}=0\leq 1-x_i$ for $i\in S(\boldsymbol x)$, whereas every
remaining $i\in I(\tilde{\boldsymbol y})$ has $x_i=0$ by the definition of
$S(\boldsymbol x)$. Hence
\[
\xi(\boldsymbol x)
\geq
\varphi\bigl(\tilde{\boldsymbol y}^{\,S(\boldsymbol x)}\bigr)
\geq
\varphi(\tilde{\boldsymbol y})-\sum_{i\in S(\boldsymbol x)}M_i(\tilde{\boldsymbol y})
=
\varphi(\tilde{\boldsymbol y})-\sum_{i\in I(\tilde{\boldsymbol y})}M_i(\tilde{\boldsymbol y})x_i,
\]
where the second inequality is the hypothesis. Since $t\geq\xi(\boldsymbol x)$, the
inequality~\eqref{F_cons1} associated with $\tilde{\boldsymbol y}$ holds, and
$\tilde{\boldsymbol y}\in Y$ was arbitrary.}

\red{Conversely, assume that $\boldsymbol x$ satisfies every inequality~\eqref{F_cons1} and
let $\hat{\boldsymbol y}$ be an optimal follower response to $\boldsymbol x$, so that
$\hat{\boldsymbol y}\in Y$, $\hat y_i\leq 1-x_i$ for every $i\in I$, and
$\varphi(\hat{\boldsymbol y})=\xi(\boldsymbol x)$; such a response exists because $Y$ is
finite and contains $\boldsymbol 0$. No facility selected in $\hat{\boldsymbol y}$ is
blocked, that is, $x_i=0$ for every $i\in I(\hat{\boldsymbol y})$. The inequality
associated with $\tilde{\boldsymbol y}=\hat{\boldsymbol y}$ therefore reduces to
$t\geq\varphi(\hat{\boldsymbol y})=\xi(\boldsymbol x)$.}
\end{myproof}

\begin{repeatproposition}

\Fabio{For each $\tilde{\boldsymbol y}\in Y$, with $Y$ defined by the upper-budget constraint~\eqref{OF_cons_l2}, taking
$M_i(\tilde{\boldsymbol y}):=\sum_{j\in J_i:\,r_j(\tilde{\boldsymbol y})=p_{ij}}d_jp_{ij}$
for every $i\in I(\tilde{\boldsymbol y})$ makes~\eqref{F_cons1} the valid family of interdiction cuts}
\begin{equation}
    t \geq 
    \varphi(\tilde{\boldsymbol y})
    -
    \sum_{i\in I(\tilde{\boldsymbol y})}
    \Fabio{\Biggl(\,\sum_{j\in J_i:\,r_j(\tilde{\boldsymbol y})=p_{ij}}d_jp_{ij}\Biggr)}
    x_i,
    \qquad \tilde{\boldsymbol y}\in Y.
\end{equation}
\end{repeatproposition}

\begin{myproof}
\Fabio{Throughout the proof we write $M_i^-(\tilde{\boldsymbol y})$ for the coefficient in the statement. For any subset $S\subseteq I(\tilde{\boldsymbol y})$, let $\tilde{\boldsymbol y}^{\,S}$ denote the surviving follower solution obtained by setting $\tilde y_i^{\,S}=0$ for $i\in S$ and $\tilde y_i^{\,S}=\tilde y_i$ otherwise. We first establish the loss bound $\varphi(\tilde{\boldsymbol y}^{\,S})\geq\varphi(\tilde{\boldsymbol y})-\sum_{i\in S}M_i^-(\tilde{\boldsymbol y})$, from which the validity of the cuts follows \red{due to Proposition~\ref{prop:validity}}.}

\red{Fix $\tilde{\boldsymbol y}\in Y$ and $S\subseteq I(\tilde{\boldsymbol y})$, and consider a customer $j\in J$. Deleting the facilities of $S$ replaces the index set $I_j(\tilde{\boldsymbol y})$ by $I_j(\tilde{\boldsymbol y})\setminus S$, so by~\eqref{OF_cons_l1} the coverage of $j$ becomes the minimum of $p_{ij}$ over a subset of the facilities over which it was previously taken. A minimum over a smaller index set cannot be smaller. Hence
\[
r_j(\tilde{\boldsymbol y}^{\,S})\geq r_j(\tilde{\boldsymbol y})
\quad\text{whenever } I_j(\tilde{\boldsymbol y})\not\subseteq S,
\qquad
r_j(\tilde{\boldsymbol y}^{\,S})=0
\quad\text{whenever } I_j(\tilde{\boldsymbol y})\subseteq S,
\]
the second case being the one in which customer $j$ is left with no selected facility at all.}

\red{Only customers of the second kind can lose coverage, and such a customer loses exactly $d_jr_j(\tilde{\boldsymbol y})$. This loss is charged to $S$ by the coefficients: if $I_j(\tilde{\boldsymbol y})\subseteq S$ is nonempty, pick $i^\ast\in\arg\min_{i\in I_j(\tilde{\boldsymbol y})}p_{ij}$, so that $i^\ast\in S$ and $r_j(\tilde{\boldsymbol y})=p_{i^\ast j}$; facility $i^\ast$ is then critical for $j$ and the term $d_jp_{i^\ast j}=d_jr_j(\tilde{\boldsymbol y})$ occurs in $M_{i^\ast}^-(\tilde{\boldsymbol y})$. If $I_j(\tilde{\boldsymbol y})$ is empty the loss is zero and nothing has to be charged. The terms charged to distinct customers are distinct summands of $\sum_{i\in S}M_i^-(\tilde{\boldsymbol y})$, because they carry different indices $j$, and all summands are nonnegative. Every customer that loses coverage is therefore charged at least once, so that
\[
\varphi(\tilde{\boldsymbol y})-\varphi(\tilde{\boldsymbol y}^{\,S})
=\sum_{j\in J}d_j\bigl(r_j(\tilde{\boldsymbol y})-r_j(\tilde{\boldsymbol y}^{\,S})\bigr)
\leq\sum_{j\in J:\,I_j(\tilde{\boldsymbol y})\subseteq S}d_jr_j(\tilde{\boldsymbol y})
\leq\sum_{i\in S}M_i^-(\tilde{\boldsymbol y}),
\]
which is the required loss bound. The first inequality drops the customers whose coverage does not decrease, and the second is the charging argument just described.}
\end{myproof}

\begin{repeatproposition}

For a fixed solution $\tilde{\boldsymbol y}$, define the dual multipliers componentwise as
\begin{equation}
\label{closedform_alpha_ec}
\tilde{\alpha}_{ij}
=
\begin{cases}
d_j,
& \text{if } j\in J_Q \text{ and } i=i(j),\\
0,
& \text{otherwise},
\end{cases}
\qquad j\in J,\ i\in I_j,
\end{equation}
and
\begin{equation}
\label{closedform_gamma_ec}
\tilde{\gamma}_j
=
\begin{cases}
0,
& \text{if } j\in J_Q,\\
d_j,
& \text{if } j\in J_N,
\end{cases}
\qquad j\in J.
\end{equation}
Then $(\tilde{\boldsymbol\alpha},\tilde{\boldsymbol\gamma})$ is an optimal solution to the full dual problem~\eqref{dualLP}.
\end{repeatproposition}

\begin{myproof}
The multipliers in~\eqref{closedform_alpha_ec}--\eqref{closedform_gamma_ec} are nonnegative by construction. For every $j\in J_Q$, equation~\eqref{closedform_alpha_ec} assigns the full dual mass $d_j$ to the minimizing index $i(j)$ while~\eqref{closedform_gamma_ec} assigns zero to $\tilde\gamma_j$; for every $j\in J_N$, equation~\eqref{closedform_gamma_ec} assigns the full mass $d_j$ to $\tilde\gamma_j$. Hence, the dual constraints~\eqref{dualLP_cons} hold at equality for every customer, and the multipliers are feasible for~\eqref{dualLP}. \blue{The degenerate case $I_j=\varnothing$ is covered as well. There is then no multiplier $\alpha_{ij}$ for customer $j$ and~\eqref{dualLP_cons} reduces to $\gamma_j\geq d_j$; moreover $\tilde N_j=0\leq 1=\tilde Q_j$ by~\eqref{MNj}, so that $j\in J_N$ and~\eqref{closedform_gamma_ec} sets $\tilde\gamma_j=d_j$, as required.}

Substituting~\eqref{closedform_alpha_ec}--\eqref{closedform_gamma_ec} into the dual objective~\eqref{dualLP_obj} gives
\[
\sum_{j\in J_Q}d_j\tilde Q_j
+
\sum_{j\in J_N}d_j\tilde N_j
=
\sum_{j\in J}d_j\min\{\tilde Q_j,\tilde N_j\}.
\]
The last expression is the optimal value of the primal projection problem~\eqref{BendersSub}. \blue{Therefore $(\tilde{\boldsymbol\alpha},\tilde{\boldsymbol\gamma})$ is dual feasible with an objective value equal to the optimal value of the primal maximization problem, and weak duality already forces it to be optimal for~\eqref{dualLP}: no dual-feasible point can have a smaller objective value without violating the weak-duality bound.}
\end{myproof}

\section{Detailed computational results}
\phantomsection
\label{app:detailed-results}
This appendix documents the runs behind the aggregate figures of Table~\ref{tab:benders-aggregate}, so that they can be reproduced instance by instance. The complete instance-level results are reported in Table~\ref{tab:main-results}. Each parameter setting contains five independently generated instances differing only in the random seed. For \red{\algname{OD+MILP}}, we report the incumbent, best bound, final optimality gap, and solution time; for \red{\algname{OD+B\&BC}}, we report the objective value and solution time. For a time-limited \red{\algname{OD+B\&BC}} run, the objective is shown as the incumbent/best-bound pair, with the final gap reported together with TL (time limit) in the time column.
\DetailedResultsTable

\section{Supplementary computational diagnostics}
\phantomsection
\label{app:additional-results}
This appendix reports five further computational experiments. We evaluate a degeneracy-based variant of the Benders separation, decompose the running time of the final algorithm into its components, measure the sensitivity of the algorithm to the target parameter and the follower budget, quantify how often several facilities are simultaneously critical for the same customer, and examine whether root-node fractional Benders separation provides a computational benefit. Together they characterize the behavior of the configuration adopted in Section~\ref{sec:DoE}.

\subsection{Dual-degeneracy experiment}
\phantomsection
\label{app:degeneracy-results}
When the Benders dual admits several optimal solutions, one may generate an additional cut from a second optimal dual vector in the hope of strengthening the Benders master problem. At a binary incumbent this happens exactly for the customers reached by a single selected facility $i(j)$, for which $\tilde Q_j=\tilde N_j$: the standard cut charges them to $\tilde\gamma_j$, the additional cut to $\tilde\alpha_{i(j)j}$, replacing the sum-based term of~\eqref{CBC} by the min-based one for all such customers at once. The additional cut is added only when it is violated by the incumbent. We evaluate this degeneracy-based variant against the standard closed-form separation on the instances solved by both. Table~\ref{tab:degeneracy-ablation-full} reports the results. The variant produces between 1.7 and 2.0 times as many cuts, but the branch-and-bound trees are of comparable size and the effect on running time depends on the instance class: it is mildly favorable for $|I|=20$ and $|I|=25$, with geometric-mean time ratios of 0.92 and 0.98, and clearly unfavorable for $|I|=35$, with a ratio of 1.70 and no win in twelve runs. Since the gain does not persist on the hardest instances, the standard separation is adopted.
\begin{table}[!htbp]
\centering
\small
\caption{Effect of dual-degeneracy cuts}
\label{tab:degeneracy-ablation-full}

\setlength{\tabcolsep}{5.0pt}
\renewcommand{\arraystretch}{0.95}

\begin{tabular}{rrrrrr}
\toprule
$|I|$
& \# common
& \shortstack{gmean time\\ratio}
& \shortstack{wins\\(Closed-form/Deg./tie)}
& \shortstack{gmean cut\\ratio}
& \shortstack{gmean max.-node\\ratio} \\
\midrule
20 & 25 & 0.92 & 6/18/1  & 1.82 & 1.03 \\
25 & 25 & 0.98 & 11/14/0 & 2.00 & 1.05 \\
35 & 12 & 1.70 & 12/0/0  & 1.74 & 0.96 \\
\bottomrule
\end{tabular}

\vspace{2pt}
\begin{minipage}{0.92\textwidth}
\footnotesize\emph{Note.}
The time, cut-count, and max.-node ratios are geometric means of the instance-wise degeneracy-based/standard ratios of the branch-and-Benders-cut algorithm over commonly solved instances. Deg.\ denotes the degeneracy-based variant.
\end{minipage}
\end{table}

\subsection{Computational time breakdown}
\phantomsection
\label{app:time-breakdown}
To locate the computational bottleneck of the final algorithm, we instrument a complete run and attribute the total time to the \red{relaxed master} solves, the separation routine, the submission of the generated cuts, the remaining solver-side effort in the follower, and the bookkeeping outside the two solves. Table~\ref{tab:time-breakdown-full} reports the breakdown on nine representative instances. Separation accounts for a substantial share of the total when the facility dimension is moderate, 29.9\% at $(|I|,|J|)=(20,20000)$ and 42.3\% at $(25,20000)$, but its share falls to 10.1\% at $(35,10000)$, where the solver-side follower effort reaches 87.7\%. The \red{relaxed master problem} and the outer bookkeeping never exceed 4.5\% and 1.0\% respectively. The measurements therefore confirm that, once the repeated auxiliary LP solves have been removed by the closed-form construction, the dominant cost is the branch-and-bound search of the follower, and that this shift occurs as the number of candidate facilities grows.
\begin{table}[!htbp]
\centering
\scriptsize
\caption{Computational time breakdown of the final algorithm}
\label{tab:time-breakdown-full}

\setlength{\tabcolsep}{4.pt}
\renewcommand{\arraystretch}{0.95}

\begin{tabular}{rrrrrrrrr}
\toprule
$|I|$ & $|J|$
& \shortstack{total\\time (s)}
& \shortstack{follower\\solves}
& \shortstack{outer master\\s (\%)}
& \shortstack{separation\\s (\%)}
& \shortstack{cut submission\\s (\%)}
& \shortstack{solver residual\\s (\%)}
& \shortstack{outer overhead\\s (\%)} \\
\midrule

20 & 20000
& 15.36 & 155
& 0.69 (4.49)
& 4.59 (29.88)
& 0.29 (1.90)
& 9.63 (62.74)
& 0.15 (0.99) \\

25 & 20000
& 49.28 & 186
& 1.07 (2.18)
& 20.85 (42.32)
& 1.36 (2.77)
& 25.79 (52.34)
& 0.20 (0.40) \\

35 & 10000
& 3250.52 & 402
& 24.21 (0.74)
& 327.64 (10.08)
& 45.18 (1.39)
& 2851.06 (87.71)
& 2.43 (0.07) \\

\bottomrule
\end{tabular}

\vspace{2pt}
\begin{minipage}{0.98\textwidth}
\footnotesize\emph{Note.}
Each row aggregates three independently generated instances. ``Separation'' measures the time spent inside the callback excluding cut materialization and submission. ``Solver residual'' is the follower solution time excluding the measured callback time and therefore includes solver-side operations that cannot be isolated reliably. ``Outer overhead'' is the residual time outside the \red{relaxed}-master and follower solves, including model updates, interdiction cut construction, and other bookkeeping operations.
\end{minipage}
\end{table}

\subsection{Sensitivity to $\beta$ and $b$}
\phantomsection
\label{app:sensitivity-results}
The baseline experiments fix the target fraction at $\beta=0.7$ and the follower budget at $b=\lfloor 0.2m\rfloor$. We vary each parameter in turn on three instances with $(|I|,|J|)=(25,10000)$, keeping the other at its baseline value. Table~\ref{tab:sensitivity-beta-b-full} reports the results. Raising $\beta$ relaxes the requirement imposed on the leader and the optimal blocker objective falls monotonically, from 10.67 at $\beta=0.6$ to 3.33 at $\beta=0.9$; the running time follows the same trend at the upper end, dropping to 3.21 seconds at $\beta=0.9$. The follower budget behaves differently. The unblocked follower value saturates at $b=7$, and so does the optimal blocker objective, which stays at 6.33 for $b\in\{7,12,17\}$; the max--min criterion means that further budget buys the follower nothing. The computational effort, however, keeps growing, from 12.96 seconds and 3,182 nodes at $b=5$ to about 100 seconds and up to 58,580 nodes at $b=17$, because the follower search space widens even though its optimal value does not change.

\begin{table}[!htbp]
\centering
\small
\caption{Detailed sensitivity results for $\beta$ and the follower budget $b$}
\label{tab:sensitivity-beta-b-full}
\setlength{\tabcolsep}{4.5pt}
\renewcommand{\arraystretch}{0.95}

\begin{tabular}{ccrrrrrr}
\toprule
Parameter
& Value
& Avg. $D^0$
& Avg. obj.
& Solved
& Time (s)
& Cuts
& Max. nodes \\
\midrule

\multirow{4}{*}{$\beta$}
& 0.6
& --
& 10.67
& 3/3
& 9.86
& 4,725
& 2,145 \\

& 0.7$^\ast$
& --
& 8.33
& 3/3
& 12.96
& 7,351
& 3,182 \\

& 0.8
& --
& 6.33
& 3/3
& 10.69
& 6,051
& 3,230 \\

& 0.9
& --
& 3.33
& 3/3
& 3.21
& 1,916
& 2,669 \\

\midrule

\multirow{4}{*}{$b$}
& 5$^\ast$
& 14965.77
& 8.33
& 3/3
& 12.96
& 7,351
& 3,182 \\

& 7
& 16981.53
& 6.33
& 3/3
& 68.50
& 26,097
& 12,353 \\

& 12
& 16981.53
& 6.33
& 3/3
& 101.35
& 28,506
& 50,309 \\

& 17
& 16981.53
& 6.33
& 3/3
& 98.00
& 24,793
& 58,580 \\

\bottomrule
\end{tabular}

\vspace{2pt}
\begin{minipage}{0.98\linewidth}
\footnotesize
\emph{Note.}
$^\ast$ denotes the baseline setting $(\beta,b)=(0.7,5)$.
Avg. $D^0$ is the arithmetic mean of the unblocked follower value over the three seeds and is reported for the budget experiment; it is unchanged across the $\beta$ experiment because $b=5$ is fixed.
Avg. obj. is the arithmetic mean of the optimal blocker objective over the three seeds.
Time, cuts, and Max. nodes are geometric means.
Max. nodes denotes the maximum number of branch-and-bound nodes observed in any individual follower solve.
\end{minipage}
\end{table}

\end{document}